\documentclass[11pt,a4paper]{article}  %
\usepackage{ifdraft}

\newcommand{\bibfile}{\jobname.bib}  %
\newcommand{\universalbib}{/home/zaikunzhang/Bureau/bibliographie/ref.bib}
\ifdraft{\IfFileExists{\universalbib}{\renewcommand{\bibfile}{\universalbib}}{}}{}
\newcommand{\iscite}{0}  %
\def\oldcite{} \let\oldcite=\cite \def\cite{\renewcommand\iscite{1}\oldcite}

\usepackage{amsmath,amsthm,amssymb,amsfonts}
\usepackage{mathtools}
\usepackage{empheq}
\usepackage{xcolor}
\usepackage[bbgreekl]{mathbbol}
\DeclareSymbolFontAlphabet{\mathbbm}{bbold}
\DeclareSymbolFontAlphabet{\mathbb}{AMSb}
\usepackage{bbm}
\usepackage{upgreek}
\usepackage{accents}
\usepackage{xspace}
\usepackage{rotating}
\usepackage{multirow,booktabs}
\usepackage[en-US]{datetime2}

\usepackage[normalem]{ulem}
\usepackage[toc,page]{appendix}

\usepackage{sectsty}
\sectionfont{\large}
\subsectionfont{\large}

\definecolor{darkblue}{rgb}{0,0.1,0.5}
\definecolor{darkgreen}{rgb}{0,0.5,0.1}
\definecolor{darkyellow}{rgb}{0.65,0.65,0.01}

\ifdraft{
\setlength{\marginparwidth}{2.42cm}
\usepackage[
tickmarkheight=3pt,
textsize=small,
backgroundcolor=blue!16,
linecolor=purple,
bordercolor=purple
]{todonotes}
}{
\newcommand{\todo}[1]{}

}

\usepackage{graphicx}
\usepackage{tikz,tikzscale,pgf}
\usetikzlibrary{arrows,arrows.meta,patterns,positioning,decorations.markings,shapes}
\usepackage{pgfplots}
\usepackage{pgfplotstable}
\usepackage[justification=centering]{caption}
\usepgfplotslibrary{fillbetween}
\pgfplotsset{compat=1.11}

\ifdraft{
\usepackage{silence}
\WarningFilter{xcolor}{Incompatible color definition on}
\WarningFilter{hyperref}{Draft mode on}
\WarningFilter{refcheck}{Unused label}
\WarningFilter{microtype}{`draft' option active}
\WarningFilter{latex}{Writing or overwriting file} %
\WarningFilter{latex}{Writing file} %
\WarningFilter{latex}{Tab has} %
\WarningFilter{latex}{Marginpar on page} %
\WarningFilter{latex}{author given} %
}{}

\ifdraft{\usepackage{refcheck}
\newcommand{\url}{\texttt}
}{
\usepackage{hyperref}
\hypersetup{colorlinks, linkcolor=darkblue, anchorcolor=darkblue, citecolor=darkblue, urlcolor=darkblue}
\usepackage{url}
}%

\usepackage{eqlist}
\usepackage{enumitem}
\setlist[itemize]{leftmargin=*}
\setlist[enumerate]{leftmargin=*,label=\normalfont{(\alph*)}}

\newcommand\overmat[2]{%
\makebox[1pt][l]{$\smash{\overbrace{\phantom{\begin{matrix}#2\end{matrix}}}^{\text{$#1$}}}$}#2}

\usepackage[section]{algorithm}
\usepackage{algpseudocode,algorithmicx}

\algrenewcommand\algorithmicrequire{\textbf{Input:}}
\algrenewcommand\algorithmicensure{\textbf{Output:}}
\algrenewcommand\alglinenumber[1]{\normalsize #1.}

\newtheorem{theorem}{Theorem}%
\newtheorem{lemma}{Lemma}%
\newtheorem{proposition}{Proposition}%
\theoremstyle{definition}
\newtheorem{remark}{Remark}%
\usepackage{xpatch}
\xpatchcmd{\proof}{\itshape}{\normalfont\proofnamefont}{}{}
\newcommand{\proofnamefont}{\bfseries}

\usepackage{microtype}
\usepackage[nobottomtitles*]{titlesec} %
\mathcode`@="8000{\catcode`\@=\active\gdef@{\mkern1mu}}

\usepackage{relsize}
\usepackage{nccmath}

\DeclareMathOperator{\tr}{tr}

\DeclareMathOperator*{\Argmin}{Argmin}

\DeclareMathOperator{\rank}{rank}
\DeclareMathOperator{\range}{range}

\DeclareMathOperator{\dist}{dist}

\newcommand{\ind}{\mathbbm{1}}

\newcommand{\QQ}{\mathcal{Q}}
\newcommand{\RR}{\mathbb{R}}
\renewcommand{\AA}{\mathcal{A}}  %
\newcommand{\BB}{\mathcal{B}}

\newcommand{\TT}{\mathcal{T}}

\newcommand{\sset}{\mathcal{S}}
\newcommand{\res}{\rho}
\newcommand{\pen}{\res}
\newcommand{\penpar}{\mu}
\newcommand{\col}{r}

\newcommand{\stf}[1]{\mathbb{S}^{#1}}
\newcommand{\sss}[1]{{\scriptscriptstyle{#1}}}

\newcommand{\fro}{{\scriptstyle{\textnormal{F}}}}
\newcommand{\trs}{{\scriptstyle{\mathsf{T}}}}
\newcommand{\signmat}{S}
\newcommand{\tsign}{{\scriptscriptstyle{S}}}
\newcommand{\posset}{\mathcal{P}}
\newcommand{\negset}{\mathcal{N}}

\newcommand{\ones}{\mathbf{1}}

\newcommand{\hp}{\circ}

\newcommand{\ie}{{i.e.}}

\xspaceaddexceptions{]\}}

\DeclareMathAlphabet{\mathsfit}{T1}{\sfdefault}{\mddefault}{\sldefault}
\SetMathAlphabet{\mathsfit}{bold}{T1}{\sfdefault}{\bfdefault}{\sldefault}

\DeclareFontFamily{U}{dutchcal}{\skewchar\font=45 }
\DeclareFontShape{U}{dutchcal}{m}{n}{<-> s*[1.0] dutchcal-r}{}
\DeclareFontShape{U}{dutchcal}{b}{n}{<-> s*[1.0] dutchcal-b}{}
\DeclareMathAlphabet{\mathlcal}{U}{dutchcal}{m}{n}
\SetMathAlphabet{\mathlcal}{bold}{U}{dutchcal}{b}{n}

\usepackage[scr=boondoxupr]{mathalfa}

\DeclareMathAlphabet{\mathpzc}{OT1}{pzc}{m}{it} %

\newcommand{\email}[1]{\texttt{\mbox{#1}}}

\title{Tight Error Bounds for the Sign-Constrained Stiefel Manifold}

\author{Xiaojun Chen~\thanks{Department of Applied Mathematics, The Hong Kong Polytechnic University,
        Hong Kong, China. Email: \email{maxjchen@polyu.edu.hk}.
    The research of this coauthor is partially supported by Hong Kong Research Grant Council project PolyU15300021.}
    \and Yifan He~\thanks{Department of Applied Mathematics, The Hong Kong Polytechnic University, Hong
        Kong, China. Email: \email{doublefan.he@connect.polyu.hk}.}
     \and Zaikun Zhang~\thanks{Department of Applied Mathematics, The Hong Kong Polytechnic University,
         Hong Kong, China. Email: \email{zaikun.zhang@polyu.edu.hk}.
             The research of this coauthor is partially supported by Hong Kong Research
        Grant Council projects PolyU15305420, PolyU15306621, and PolyU15308623.}
        }

\begin{document}

\maketitle

\begin{abstract}
The sign-constrained Stiefel manifold in~$\RR^{n\times \col}$ is
a segment of
the Stiefel manifold with fixed signs~(nonnegative or nonpositive) for some columns of the matrices.
It includes the nonnegative Stiefel manifold as a special case.
We present global and local error bounds that provide an inequality with easily computable
residual functions and explicit coefficients to bound the distance from matrices
in~$\RR^{n\times \col}$ to the sign-constrained Stiefel manifold.
Moreover,
we show that the error bounds cannot be improved except for the multiplicative constants under some mild conditions, which explains why two
square-root terms are necessary in the bounds when~$1<\col<n$ and why the~$\ell_1$ norm can be used in
the bounds when~$\col=n$ or~$\col=1$ for the sign constraints and orthogonality, respectively.
The error bounds are applied to derive exact penalty methods for minimizing a Lipschitz continuous function with
orthogonality and sign constraints.
\end{abstract}



\section{Introduction}\label{sec:inro}

Let~$n$ and~$\col$ be two integers such that~$1\le \col \le n$, and
$\stf{n,\col}\coloneqq\{ X \in \RR^{n\times \col} \mathrel{:} X^\trs X = I_\col \}$ be the
Stiefel manifold, where~$I_\col$ is the~$\col \times \col$ identity matrix.
Given two disjoint subsets~$\posset$ and~$\negset$ of~$\{j\mathrel{:} 1\le j \le \col\}$,
denote
\begin{equation*}
    \RR^{n\times \col}_{\tsign}
    \;\coloneqq\; \left\{ X \in \RR^{n\times \col} \mathrel{:} X_{ij} \ge 0 \text{ for } j \in \posset
        \text{ and }
    X_{ij} \le 0 \text{ for } j \in \negset, ~ 1\le i \le n \right\},
\end{equation*}
which is a subset of $\RR^{n\times \col}$ with column-wise nonnegative or nonpositive constraints on some columns.

In this paper, we consider the \emph{sign-constrained Stiefel manifold}
defined as
\begin{equation}
    \nonumber
    \stf{n,\col}_{\tsign} \;\coloneqq\; \stf{n,\col} \!\cap \RR^{n\times \col}_{\tsign}.
\end{equation}
When~$\posset = \{j\mathrel{:} 1\le j \le \col\}$, $\RR^{n\times\col}_\tsign$ reduces to the
nonnegative orthant~$\RR^{n\times \col}_+$, and~$\stf{n,\col}_\tsign$ reduces to the
\emph{nonnegative Stiefel manifold}~$\stf{n,\col}_+ \coloneqq \{X\in \stf{n,\col}\mathrel{:} X \ge 0\}$.

If we define~$\signmat\in\RR^{n\times \col}$ as the matrix with
\begin{equation}
    \nonumber
    S_{i,j} \;=\;
    \begin{cases}
        1, & \text{if } j \in \posset,\\
        -1, & \text{if } j \in \negset,\\
        0, & \text{otherwise},\\
    \end{cases}
    ~\;~ 1\le i \le n,
\end{equation}
then~$\stf{n,\col}_{\tsign}$~\!can be formulated as
\begin{equation}
    \nonumber
    \stf{n,\col}_{\tsign} \;=\; \{X\in \RR^{n\times\col} \mathrel{:} \signmat \hp X \ge 0, ~X^\trs
    X = I_\col\},
\end{equation}
where~$\hp$ signifies the Hadamard product.
We will investigate error bounds
  \begin{eqnarray}
  \label{eq:error1}
  & \dist(X,\, \stf{n,\col}_\tsign) \;\le\;  \nu \|(\signmat \hp X)_{-}\|_\fro^q & \quad \text{for} \quad X \in\stf{n,\col},\\
  \label{eq:error2}
  &\dist(X,\, \stf{n,\col}_\tsign) \;\le\; \nu\|X^\trs X-I_\col\|_\fro^q& \quad \text{for} \quad
       X \in \RR^{n\times \col}_\tsign,\\
  \label{eq:error3}
  & \dist(X,\, \stf{n,\col}_\tsign)  \;\le\; \nu(\|(\signmat\hp X)_{-}\|_\fro^{q} +
  \|X^\trs X-I_\col\|_\fro^{q})
      &\quad \text{for} \quad X \in \RR^{n\times \col},
   \end{eqnarray}
where~$\nu$ and~$q$ are positive constants,
and~$Y_{-}\coloneqq \max\{-Y, \, 0\}$ stands for the entry-wise nonnegative part of~$-Y$ for any
matrix~$Y$.
The bounds~\mbox{\eqref{eq:error1}--\eqref{eq:error3}}
are global error bounds for~$\stf{n,\col}_\tsign$ relative to~$\stf{n,\col}$,
$\RR^{n\times \col}_\tsign$, and~$\RR^{n\times \col}$, respectively, with the first two being
special cases of the last one.

According to the error bound of Luo-Pang presented in~\cite[Theorem~2.2]{Luo_Pang_1994},
there exist~$\nu>0$ and~$q>0$ such that the
inequalities
in~\mbox{\eqref{eq:error1}--\eqref{eq:error3}} hold for all~$X$ in a compact subset of~$\RR^{n\times \col}$.
Moreover, due to the error bound for polynomial systems  given
in~\cite[Corollary 3.8]{Li_Mordukhovich_Pham_2015},
for all~$X$ in a compact subset of~$\RR^{n\times \col}$, there exists a~$\nu$ such
that the inequalities in~(\ref{eq:error1})--(\ref{eq:error3}) hold with a
dimension-dependent value of~$q$ that is less than~$6^{-2nr}$.
However, to the best of our knowledge, the explicit value of~$\nu$ and the value of~$q$ that is
independent of the dimension in~(\ref{eq:error1})--(\ref{eq:error3}) are still unknown
even in the special case of~$\stf{n,\col}_\tsign\!\!= \stf{n,\col}_+$,
and it is also unknown whether the error bounds hold in an unbounded set.

Being a fundamental concept in optimization, error bound plays a crucial role in
both theory and methods for solving systems of equations and optimization problems~\cite{Luo_Pang_1994,Pang_1997}.
One of its applications is to develop penalty methods for constrained optimization problems.
Let~$F\mathrel{:} \RR^{n\times \col} \to \RR$ be a continuous function.
The minimization problem
\begin{equation}
    \label{eq:optnoc}
    \min \left\{F(X) \mathrel{:} X \in \stf{n,\col}_\tsign \right\}
\end{equation}
can be found in a wide range of optimization models in data science,
including  nonnegative principal component
analysis~\cite{Lui-Boumal,Zass_Shashua_2006},
nonnegative Laplacian embedding~\cite{Luo_Ding_Huang_Li_2009}, discriminative nonnegative spectral
clustering~\cite{Yang_etal_2012},
orthogonal nonnegative  matrix factorization~\cite{Pompili_Gillis_Absil_Glineu\col_2014,Yang_Oja_2010},
and some K-indicators models for data clustering~\cite{Chen_etal_2019,WCCP2021}.

Even in the special case of~$\stf{n,\col}_\tsign\!\!= \stf{n,\col}_+$, the constraints of
problem~\eqref{eq:optnoc} are challenging to handle due to their combinatorial nature
(note that, for example,~$\stf{n,n}_+$ equals the set of all permutation matrices on~$\RR^n$).~To
deal with these difficult constraints, the penalty problems
 \begin{align}
    \label{penalty1}
    \min & \left\{F(X) + \penpar \|(S\hp X)_{-}\|_\fro^q \mathrel{:} X \in \stf{n,\col} \right\},\\
    \label{penalty2}
    \min & \left\{F(X) + \penpar\|X^\trs X-I_\col\|_\fro^q \mathrel{:} X \in \RR^{n\times \col}_\tsign\right\},\\
    \label{penalty3}
    \min  & \left\{F(X) + \penpar(\|(S\hp X)_{-}\|_\fro^{q}+ \|X^\trs X-I_\col\|_\fro^{q}) \mathrel{:} X \in \RR^{n\times \col}\right\},
\end{align}
have been widely used for solving~\eqref{eq:optnoc} with~$\stf{n,\col}_\tsign\!\!=\stf{n,\col}_+$,
where~$\penpar$ is the penalty parameter.
See for example~\mbox{\cite{Ahookhosh_etal_2021,Qian_Pan_Xiao_2023,Yang_etal_2012,Zass_Shashua_2006}}
and the references therein.
However, the exactness of
problems~\mbox{(\ref{penalty1})--(\ref{penalty3})} regarding global minimizers and local minimizers of
problem~(\ref{eq:optnoc}) is not well understood.

The main contribution of this paper is to establish
the error bounds~\mbox{(\ref{eq:error1})--(\ref{eq:error3})} with~$\nu = 15\col^{\frac{3}{4}}$
and~$q={1}/{2}$ without any additional restriction on~$X$.
Moreover, we demonstrate that the error bounds cannot hold for~$q>{1} /{2}$ under mild conditions
when~$1 < \col < n$ and~$\stf{n,\col}_\tsign=\stf{n,\col}_+$.
In addition, we show that the error bounds~(\ref{eq:error1})--(\ref{eq:error3}) hold
with~$q=1$ and~$\nu = 7\sqrt{\col}$ when~$|\posset|+|\negset| = 1$, and hold with~$q = 1$ and~$\nu = 9n$ when~$|\posset|+|\negset|=n$,
but they cannot hold with~$q>1$.
As an application of error bounds~\mbox{(\ref{eq:error1})--(\ref{eq:error3})} with~$\nu = 15\col^{\frac{3}{4}}$
and~$q={1}/{2}$, we show the exactness of the
penalty problems~(\ref{penalty1}) and~(\ref{penalty2}) under the assumption that~$F$ is Lipschitz
continuous, taking~$\stf{n,\col}_\tsign = \stf{n, \col}_+$ as an example. Moreover, we show the
existence of Lipschitz continuous functions such that  penalty problems~(\ref{penalty1})
and~(\ref{penalty2}) with~$q>1/2$ is not an exact penalty for global and local minimizers of the corresponding constrained problems.

Very recently, our error bound and matrix inequalities have been used to study constant modulus
optimization and optimal orthogonal channel selection~\cite{Channel,Anthony1, Anthony2}, which have a wide variety of applications in signal processing, communications, and data science.

The rest of the paper is organized as follows.
In Section~\ref{sec:lemmas}, we introduce some notations and preliminaries.
Section~\ref{sec:errbd} derives the error bounds~(\ref{eq:error1})--(\ref{eq:error3}) in the special
case of~$\stf{n,\col}_\tsign\!= \stf{n,\col}_+$.
Section~\ref{sec:serrbd} extends these bounds to the general case.
Section~\ref{sec:exactpen} investigates the exactness of the penalty problems~\eqref{penalty1}-\eqref{penalty3} using the new error bounds. Section~\ref{sec:test} considers applications of the theoretical results established in this paper.
We conclude the paper in Section~\ref{sec:conclusion}.%

\section{Notation and preliminaries}
\label{sec:lemmas}

For any matrix~$X\in\RR^{n\times \col}, X_+\coloneqq \max\{X, \, 0\}= X + X_{-}$ is
the projection of~$X$ onto~$\RR^{n\times \col}_+$. In addition,
the singular value vector of~$X$ is denoted by~$\sigma(X)\in\RR^{\col}$,
the entries of which are in the descent order.
Meanwhile,~$\Sigma(X)\in\RR^{n\times \col}$ is the matrix such
that~$X = U\Sigma(X)V^\trs$ is the singular value decomposition of~$X$, the diagonal
of~$\Sigma(X)$ being~$\sigma(X)$.
We use~$\ones$ to denote the vector with all entries being one, and its dimension will be clear from
the context.

Unless otherwise specified,~$\|\cdot\|$ stands for a general vector norm.
For any constant~$p \in [1, +\infty)$, we use~$\|\cdot\|_p$ to represent either
the~$\ell_p$-norm of vectors or the operator norm induced by this vector norm for matrices.
In addition, we use~$\|\cdot\|_{\ell_p}$ to denote the
\emph{entry-wise}~$\ell_p$-norm of a matrix, namely the~$\ell_p$-norm of the vector that contains
all the entries of the matrix. Note
that~$\|\cdot\|_{\ell_2}$ is the Frobenius norm, which is also denoted by~$\|\cdot\|_\fro$.
When~$\RR^{n\times \col}$ is equipped with the Frobenius norm,
we use~$\BB(X,\delta)$ to represent the open ball in~$\RR^{n\times\col}$
centered at a point~$X \in \RR^{n\times \col}$ with a radius~$\delta>0$,
and~$\dist(X,\, \TT)$ to denote the distance from a point~$X\in \RR^{n\times \col}$
to a set~$\TT\subset\RR^{n\times\col}$.
Finally, given an minimization problem, we use~$\Argmin$ to denote the set of global minimizers.

Lemma~\ref{lem:Mirsky} is fundamental for the analysis of distances between matrices.
This lemma is stated for unitarily invariant norms~(see~\cite[Section~3.5]{Horn_Johnson_2008}
for this concept), although we are most interested in the case with the Frobenius norm.

\begin{lemma}
    [Mirsky]
    \label{lem:Mirsky}
    For any matrices~$X \in \RR^{n\times \col}$ and~$Y\in\RR^{n\times \col}$, we have
    \begin{equation}
     \label{eq:Mirsky}
    \|\Sigma(X) - \Sigma(Y)\| \;\le\; \|X-Y\|
    \end{equation}
    for any unitarily invariant norm~$\|\cdot\|$
    on~$\RR^{n\times \col}$.
    When~$\|\cdot\|$ is the Frobenius norm, the equality holds in~\eqref{eq:Mirsky}
    if and only if there exist orthogonal matrices~$U\in\RR^{n\times n}$ and~$V \in \RR^{\col \times
    \col}$ such that~$X = U\Sigma(X) V^\trs$ and~$Y = U\Sigma(Y)V^\trs$.
\end{lemma}

The square case~(\ie,~$n = \col$) of inequality~\eqref{eq:Mirsky} is due to
Mirsky~\cite[Theorem~5]{Mirsky_1960}, and the general case can be found in~\cite[Theorem
7.4.9.1]{Horn_Johnson_2012}.
A direct corollary of Lemma~\ref{lem:Mirsky} is the following
Hoffman-Wielandt~\cite{Hoffman_Wielandt_2003} type bound for singular values, which is equivalent to
the von Neumann trace inequality~\cite[Theorem~I]{vonNeumann_1937}~(see also~\cite[Theorem~2.1]{Lewis_1995}).
\begin{lemma}[von Neumann]
    \label{lem:vonNeumann}
    For any matrices~$X \in \RR^{n\times \col}$ and~$Y\in\RR^{n\times \col}$,~we~have%
    \begin{equation}
        \nonumber
    \|\sigma(X) - \sigma(Y)\|_2  \;\le\; \|X-Y\|_\fro,
    \end{equation}
    and equivalently,~$\tr(X^\trs Y) \le \sigma(X)^\trs \sigma(Y)$.
\end{lemma}

The following lemma is another consequence of Lemma~\ref{lem:Mirsky}.
For this result, recall that each matrix~$X\in\RR^{n\times \col}$ has a polar decomposition in the
form of~$X = UP$, where~$U$ belongs to~$\stf{n,\col}$ and~$P=(X^\trs X)^{\frac{1}{2}}$,
with~$U$ being called a unitary polar factor of~$X$.
The square case of this lemma is due to Fan and Hoffman~\cite[Theorem~1]{Fan_Hoffman_1955}.
For the general case, see~\mbox{\cite[Theorem~8.4]{Higham_2008}},
which details a proof based on Lemma~\ref{lem:Mirsky}.%
\begin{lemma}[Fan-Hoffman]\label{lem:FanHoffman}%
    If~$U\in\RR^{n\times \col}$ is a unitary polar factor of a matrix~$X \in\RR^{n\times \col}$, then%
    \begin{equation}
        \nonumber
        \|X - U\| \;=\; \min\{\|X-V\|\mathrel{:} V\in\stf{n,\col}\}
    \end{equation}
    for any unitarily invariant norm~$\|\cdot\|$ on~$\RR^{n\times \col}$.
\end{lemma}

Lemma~\ref{lem:distm} collects a few basic facts on the distance from a matrix in~$\RR^{n\times \col}$ to~$\stf{n,\col}$.
\begin{lemma}
    \label{lem:distm}
    For any matrix~$X\in\RR^{n\times \col}$, we have
    \begin{equation}
        \nonumber
    \dist(X,\,\stf{n, \col}) \;=\;  \|\sigma(X) - \ones\|_2 \;\le\;
    \min\Big\{\|X^\trs X - I_\col\|_\fro, \; \col^{\frac{1}{4}}\|X^\trs X - I_\col\|_\fro^{\frac{1}{2}}\Big\}.
    \end{equation}
In addition,~$\|X^\trs X - I_\col\|_\fro \le (\|X\|_2+1)\|\sigma(X) - \ones\|_2$.
\end{lemma}
\begin{proof}
        Let~$U\in\stf{n,\col}$ be a unitary polar factor of~$X$. By Lemma~\ref{lem:FanHoffman},
    \begin{equation}
        \nonumber
    \dist(X,\, \stf{n,\col})
    = \|X - U\|_\fro
    = \|U^\trs (X - U)\|_\fro
    = \|(X^\trs X)^{\frac{1}{2}} - I_\col\|_\fro
    = \|\sigma(X) - \ones\|_2.
    \end{equation}
        The entry-wise inequalities~$|\sigma(X) -\ones| \le |\sigma(X)^2 -\ones| \le (\|X\|_2+1)|\sigma(X) -\ones|$ imply
    \begin{equation}
    \label{eq:sigminus1}
    \|\sigma(X) -\ones\|_2 \;\le\; \|\sigma(X)^2 - \ones\|_2
    \;\le\; (\|X\|_2+1)\|\sigma(X) -\ones\|_2.
    \end{equation}
    Noting that~$\|\sigma(X)^2 -\ones\|_2 = \|X^\trs X-I_\col\|_\fro$, we obtain from~\eqref{eq:sigminus1}
    that
    \begin{equation*}
        \|\sigma(X) -\ones\|_2 \le\|X^\trs X-I_\col\|_\fro \le(\|X\|_2+1)\|\sigma(X) - \ones\|_2.
    \end{equation*}
    Finally, since~$|\sigma(X) -\ones|^2 \le |\sigma(X)^2 -\ones|$, we have
    \begin{equation}
    \nonumber
        \|\sigma(X) -\ones\|_2^2 \;\le\; \|\sigma(X)^2 - \ones\|_1
    \;\le\; \sqrt{\col} \|\sigma(X)^2 - \ones\|_2 \;=\; \sqrt{\col} \|X^\trs X -I_\col\|_\fro.
    \end{equation}
    The proof is complete.
\end{proof}

By Lemmas~\ref{lem:FanHoffman} and~\ref{lem:distm},~$\dist(X,\, \stf{n,\col}_+) = \|\sigma(X) - \ones\|_2$
if~$X$ has a nonnegative unitary polar factor.
It is the case in the following lemma, where this factor is~$X(X^\trs X)^{-\frac{1}{2}}$.%
\begin{lemma}
    \label{lem:diag}
    For a matrix~$X\in\RR^{n\times \col}_+$, if~$X^\trs X$ is nonsingular and diagonal,~then%
\begin{equation}
    \nonumber
    \dist(X,\, \stf{n,\col}_+) \;=\; \|\sigma(X) - \ones\|_2.
\end{equation}
\end{lemma}

Lemma~\ref{lem:row} is an elementary property of~$\stf{n,\col}_+$. We omit the proof.
\begin{lemma}
    \label{lem:row}
    For a matrix~$X\in\stf{n,\col}_+$, each row of~$X$ has at most one nonzero entry.
\end{lemma}

\section{Error bounds for the nonnegative Stiefel manifold}
\label{sec:errbd}

This section will establish the error bounds~\eqref{eq:error1}--\eqref{eq:error3}
 in the special case of~$\stf{n,\col}_\tsign = \stf{n,\col}_+$, where~$S\hp X$ reduces to~$X$.
Subsection~\ref{ssec:errbd1n} demonstrates~\eqref{eq:error1}--\eqref{eq:error3}
with~$q = 1$ when~$\col=1$ or~$\col=n$, and points out that
they cannot hold with~$q>1$ regardless of~$\col\in\{1,\dots, n\}$.
In Subsection~\ref{ssec:errbdr}, we derive the bounds~\eqref{eq:error1}--\eqref{eq:error3}
with~$q={1}/{2}$ for~$1\le \col \le n$, and
Subsection~\ref{ssec:tight} elaborates on the tightness of these bounds when~$1<\col<n$.
As an application of our results, we briefly discuss the linear regularity of~$\RR^{n\times \col}_+$
and~$\stf{n,\col}$ in Subsection~\ref{ssec:linreg}.

General discussions on error bounds can be found in~\cite[Section~6.1]{Facchinei_Pang_2003a}.
Here we focus on error bounds for~$\stf{n,\col}_+$ defined by two special functions
\begin{align*}
    \res_1(X) \coloneqq\; &\|X_{-}\|_\fro^{q_1} + \|\sigma(X)-\ones\|_2^{q_2},\\
    \res_2(X) \coloneqq\; & \|X_{-}\|_\fro^{q_1} + \|X^\trs X-I_\col\|_\fro^{q_2},
\end{align*}
where~$q_1$ and~$q_2$ are positive constants.
These functions are residual functions for~$\stf{n,\col}_+$ relative to~$\RR^{n\times \col}$,
namely nonnegative-valued functions on~$\RR^{n\times \col}$ whose zeros coincide with the
elements of~$\stf{n,\col}_+$.
The residual function~$\res_2$ is easily computable and it reduces to the one in~\eqref{eq:error3} when~$q_1 = q_2 = q$.

We say that~$\res_1$ defines a local error bound for~$\stf{n,\col}_+$ relative
to~$\RR^{n\times\col}$ if there exist positive constants~$\epsilon$ and~$\nu$ such that
\begin{equation}
    \label{eq:errbd}
    \dist(X,\, \stf{n,\,\col}_+) \;\le\; \nu \res_1(X)
\end{equation}
for all~$X\in \RR^{n\times \col}$ satisfying~$\|X_{-}\|_\fro+\|X^\trs X-I_{\col}\|_\fro\le \epsilon$,
and we say it defines a global error bound for~$\stf{n,\col}_+$ relative to~$\RR^{n\times\col}$
if~\eqref{eq:errbd} holds for all~$X\in\RR^{n\times \col}$.
Likewise, we can use~$\res_1$ to define error bounds for~$\stf{n,\col}_+$ relative to
any set~$\sset\subset\RR^{n\times\col}$ that contains~$\stf{n, \col}_+$, for
example,~$\sset = \RR^{n\times\col}_+$, in which case~$\res_1$ reduces to its second term.
Similar things can be said about~$\res_2$.
Theorems~\ref{th:power1n} and~\ref{th:power}
will specify the precise range of~$q_1$ and~$q_2$ so that~$\res_1$ and~$\res_2$ define local or
global error bounds for~$\stf{n,\col}_+$ relative to~$\RR^{n\times \col}$.

\subsection{Tight error bounds with~$\col =  1$ or~$\col = n$}
\label{ssec:errbd1n}

In this subsection, we show that the error bounds~(\ref{eq:error1})--(\ref{eq:error3}) hold for~$q=1$
when~$\col=1$ or~$\col=n$. Moreover, we explain why
bounds~(\ref{eq:error1})--(\ref{eq:error3}) cannot hold for~$q>1$ in general.

The bound for~$\col = 1$ is easy to establish due to the simple fact that
\begin{equation}
     \label{eq:keyeq1}
        \dist(x,\,\stf{n,1}_+) \;=\;  \dist(x,\,\stf{n,1}) \;=\;  \big|\|x\|_2 - 1\big|
        \quad \text{ for all } \quad x\in\RR^{n}_+.
\end{equation}
Indeed, when~$x = 0$, this is trivial; when~$x\neq 0$, equality~\eqref{eq:keyeq1} is true
because the projection of~$x$ onto~$\stf{n,1}_+$ equals its projection onto~$\stf{n,1}$, which
is~$x/\|x\|_2 \ge 0$.
\begin{theorem}
    \label{th:gerrbdr1}
    For any vector~$x\in\RR^{n}$,
\[
    \dist(x,\, \stf{n,1}_+)  \;\le\;  2\|x_{-}\|_2 + \big|\|x\|_2 - 1\big|.
\]
\end{theorem}
\begin{proof}
    As observed above,~$\dist(x_+,\,\stf{n,1}_+) = \big| \|x_+\|_2-1\big |$.
    Meanwhile,
    \begin{equation}
        \nonumber
        \begin{split}
        \big|\|x_{+}\|_2-1\big| - \big|\|x\|_2-1\big|
        \;\le\; \big|\|x_{+}\|_2 - \|x\|_2\big|
        \;\le\; \|x_{+}- x\|_2
        \;=\; \|x_{-}\|_2.
        \end{split}
    \end{equation}
    Thus~$\dist(x,\, \stf{n,1}_+)  \le\|x_{-}\|_2+ \dist(x_+,\, \stf{n,1}_+) \le  2\|x_{-}\|_2 + \big|\|x\|_2 - 1\big|$.
\end{proof}

To establish the error bounds for~$\col = n$, we first prove Proposition~\ref{prop:lerrbdnr+}, which
is essentially a weakened version of the observation~\eqref{eq:keyeq1} in the current situation.
Note that the matrix~$Y$ defined in the proof below is indeed the rounding matrix proposed
in~\cite[Procedure~3.1]{Jiang_Meng_Wen_Chen_2023}.
\begin{proposition}
    \label{prop:lerrbdnr+}
    For any matrix~$X\in\RR^{n\times n}_+$, if~$\|\sigma(X) - \ones\|_2 < 1/(4\sqrt{n})$, then
    \begin{equation}
    \label{eq:lerrbdnr+}
        \dist(X,\,\stf{n,n}_+)  \;\le\;7\sqrt{n} \|\sigma(X) - \ones\|_2.
    \end{equation}
\end{proposition}
\begin{proof}
    For each~$i \in \{1, \dots, n\}$, take the smallest integer~$l_i \in \{1,\dots, \col\}$ so that
    \begin{equation}
        \nonumber
        X_{i,l_i} \;=\; \max\, \{X_{i,j} \mathrel{:}  j = 1, \dots, \col\}.
    \end{equation}
    Consider the matrix~$Y \in \RR^{n,\col}_+$ defined by
    \begin{equation}
        \label{eq:y}
        Y_{i, j} \;=\;
        \begin{cases}
            X_{i,l_i} & \text{if } j = l_i,\\
            0 & \text{otherwise}.
        \end{cases}
    \end{equation}
    We will demonstrate~\eqref{eq:lerrbdnr+} by establishing bounds for~$\|X-Y\|_\fro$ and~$\dist(Y,\,\stf{n,n}_+)$.

    Consider~$\|X-Y\|_\fro$ first.
    Due to the fact that~$\|\sigma(X) - \ones\|_2 <1/(4\sqrt{n})$, all the~$n$ singular values
    of~$X$ are at least~$3/4$.
    Since~$X\ge 0$ and~$X_{i,l_i} = \max\{X_{i,j} \mathrel{:} j = 1, \dots, n\}$, we have
    \begin{equation*}
    X_{i, l_i}  \;\ge\; \frac{1}{\sqrt{n}}\left( XX^\trs \right)_{i,i}^{\frac{1}{2}}
    \;\ge\; \frac{3}{4\sqrt{n}}
        \quad\text{for each} \quad i \in\{1, \dots, n\}.
    \end{equation*}
    Fix an integer~$j\in\{1, \dots, \col\}$.
    For each~$l \in \{1, \dots, \col\}$, define
    \begin{equation}
        \label{eq:1lj}
        \nonumber
        \ind(j \neq l) \;=\;
        \ind(l \neq j) \;=\;
        \begin{cases}
            1 & \text{ if } l \neq j, \\
            0 & \text{ if } l = j.
        \end{cases}
    \end{equation}
    With~$x^j$ and~$y^j$ denoting the~$j$th columns of~$X$ and~$Y$, respectively, we have
    \begin{equation*}
        \begin{split}
            \frac{9}{16n}\|x^j - y^j\|_2^2
            \;=\;&  \frac{9}{16n} \sum_{i = 1}^{n} X_{i, j}^2\ind(j \neq l_i) \\
            \;\le\;&  \sum_{i = 1}^{n} X_{i,l_i}^2 X_{i, j}^2\ind(l_i \neq j) \\
            \;\le\;&
            \sum_{l=1}^n\sum_{i = 1}^{n} X_{i,l}^2 X_{i,j}^2 \ind(l \neq j) \\
            \;\le\;&  \sum_{l=1}^n\left(\sum_{i = 1}^{n} X_{i,l} X_{i,j} \right)^2 \ind(l \neq j) \\
            \;=\;&  \sum_{l=1}^n\left( X^\trs X -I_n \right)^2_{l,j} \ind(l \neq j).
        \end{split}
    \end{equation*}
       Hence
    \begin{equation*}
    \|X - Y\|_\fro \;\le\;
    \frac{4}{3}\sqrt{n}\|X^\trs X -I_n\|_\fro.
    \end{equation*}
     By Lemma~\ref{lem:distm} and the fact that~$\|X\|_2 \le 1 + \|\sigma(X) -\ones\|_2 \le 5/4$,
     we have further
    \begin{equation}
        \label{eq:xdysigma}
        \|X-Y\|_\fro \;\le\; \frac{4}{3} \sqrt{n}  (\|X\|_2+ 1)\|\sigma(X) -\ones\|_2
        \;\le\; 3\sqrt{n}\|\sigma(X) -\ones\|_2.
    \end{equation}

    Now we estimate~$\dist(Y,\,\stf{n,n}_+)$. According to inequality~\eqref{eq:xdysigma} and
    Lemma~\ref{lem:vonNeumann},
    \begin{equation}
        \label{eq:sydone}
        \|\sigma(Y) - \ones\|_2 \;\le\; \|X-Y\|_\fro + \|\sigma(X) - \ones\|_2 \;\le\; 4\sqrt{n}
        \|\sigma(X) -\ones\|_2.
    \end{equation}
    Since~$\|\sigma(X)-\ones\|_2 < 1/(4\sqrt{n})$, we have~$\|\sigma(Y) - \ones\|_2 < 1$, which
    implies that~$Y^\trs Y$ is nonsingular.
    Since~$Y$ has at most one nonzero entry in each row, it is clear that~$Y^\trs Y$ is diagonal.
    Thus we can invoke Lemma~\ref{lem:diag} and obtain
    \begin{equation}
        \nonumber
        \dist(Y,\,\stf{n,n}_+) \;=\;\|\sigma(Y) -\ones\|_2.
    \end{equation}
    Therefore, combining inequalities~\eqref{eq:xdysigma} and~\eqref{eq:sydone},
    we conclude that~\eqref{eq:lerrbdnr+} is true.
\end{proof}

Theorem~\ref{th:errbdnr} presents global and local error bounds for~$\stf{n,n}_+$ relative to~$\RR^{n\times n}$.
\begin{theorem}
    \label{th:errbdnr}%
        For any matrix~$X\!\!\in\RR^{n\times n}$, we have
        \begin{equation}
            \label{eq:gerrbdnr}
            \dist(X,\,\stf{n,n}_+) \;\le\; 9n\left( \|X_{-}\|_\fro + \|\sigma(X) -\ones\|_2 \right).
        \end{equation}
        Moreover, if~$\|X_{-}\|_\fro + \|\sigma(X) -\ones\|_2 < 1/(4\sqrt{n})$, then
        \begin{equation}
            \label{eq:lerrbdnr}
            \dist(X,\,\stf{n,n}_+) \;\le\; 8\sqrt{n}\left( \|X_{-}\|_\fro + \|\sigma(X) -\ones\|_2 \right).
        \end{equation}
\end{theorem}
\begin{proof}
    We first prove~\eqref{eq:lerrbdnr}, assuming that~$\|X_{-}\|_\fro + \|\sigma(X) -\ones\|_2 < 1/(4\sqrt{n})$.
    By Lemma~\ref{lem:vonNeumann}, this assumption ensures~$\|\sigma(X_+) -\ones\|_2< 1/(4\sqrt{n})$.
    Thus Proposition~\ref{prop:lerrbdnr+} renders%
    \begin{equation}
        \nonumber
        \begin{split}
        \dist(X_+, \, \stf{n,n}_+)
        \;\le\;  7\sqrt{n}\|\sigma(X_+)-\ones\|_2
        \;\le\;  7\sqrt{n}\left(\|X_{-}\|_\fro+\|\sigma(X)-\ones\|_2 \right),
        \end{split}
    \end{equation}
    which justifies inequality~\eqref{eq:lerrbdnr} since~$\dist(X, \, \stf{n,n}_+) \le \|X_{-} \|_\fro + \dist(X_+,\,\stf{n,n}_+)$.

    Now we consider inequality~\eqref{eq:gerrbdnr}. If~$\|X_{-}\|_\fro + \|\sigma(X) -\ones\|_2 < 1/(4\sqrt{n})$,
    then~\eqref{eq:gerrbdnr} holds due to~\eqref{eq:lerrbdnr}.
    When~$\|X_{-}\|_\fro + \|\sigma(X) -\ones\|_2 \ge 1/(4\sqrt{n})$,
    inequality~\eqref{eq:gerrbdnr} is justified~by%
    \begin{equation}
        \nonumber
        \begin{split}
            \dist(X,\,\stf{n,n}_+)
        \;\le\; & \dist(X,\,\stf{n,n}) + 2\sqrt{n} \\
        \;\le\; & \|\sigma(X)-\ones\|_2 + 8n(\|X_{-}\|_\fro + \|\sigma(X) -\ones\|_2) \\
        \;\le\; & 9n\left( \|X_{-}\|_\fro + \|\sigma(X) -\ones\|_2 \right),
        \end{split}
    \end{equation}
    where the first inequality holds because the diameter of~$\stf{n,n}$ is~$2\sqrt{n}$.
\end{proof}

\begin{remark}\label{remark1}
Since~$\big|\|x\|_2- 1\big|\le \big|\|x\|^2_2- 1\big|$
and~$\|\sigma(X) - \ones\|_2\le\|X^\trs X-I_n\|_{\fro}$,
Theorems~\ref{th:gerrbdr1} and~\ref{th:errbdnr}
imply the error bounds~(\ref{eq:error1})--(\ref{eq:error3}) with~$q=1$ for~$\col\in\{1,n\}$.
These bounds cannot be improved except for the multiplicative constants. Indeed,
for any matrix~$X\in \RR^{n\times \col}$ with~$\col\in \{1, \dots, n\}$ and~{$\|X\|_2 \le 1$},
we have
\begin{equation}
\label{eq:lowerdist}
\begin{split}
\dist(X,\, \stf{n,\col}_+)
\;\ge\; & \max\left\{\dist(X,\, \RR^{n\times \col}_+),\; \dist(X,\, \stf{n,\col})\right\} \\
\;\ge\; & \frac{1}{2}\left[\dist(X,\, \RR^{n\times \col}_+) + \dist(X,\, \stf{n,\col})\right] \\
\;=\; & \frac{1}{2}\left(\|X_{-}\|_\fro + \|\sigma(X) -\ones\|_2\right) \\
\;\ge\; & \frac{1}{2}\|X_{-}\|_\fro + \frac{1}{4}\|X^\trs X -I_\col\|_\fro,
\end{split}
\end{equation}
where the last two lines apply Lemma~\ref{lem:distm}.
This also implies that the bounds~\mbox{(\ref{eq:error1})--(\ref{eq:error3})}
cannot hold for any~$\col\in\{1, \dots,n\}$ with~$q>1$.
\end{remark}

Theorem~\ref{th:power1n} is an extension of Theorems~\ref{th:gerrbdr1} and~\ref{th:errbdnr}.
It specifies the possible exponents of~$\|X_{-}\|_\fro$ and~$\|\sigma(X) -\ones\|_2$
or~$\|X^\trs X -I_\col\|_\fro$ in local and global error bounds for~$\stf{n,\col}_+$ relative
to~$\RR^{n\times\col}$ for~$\col \in \{1, n\}$.
As we will see from~\ref{it:fro1n} of this theorem and its proof,
when~$\col = 1$ or~$\col = n$, the error bound~\eqref{eq:error1} can hold if and only if~$q \le 1$,
whereas~\eqref{eq:error2} and~\eqref{eq:error3} can hold if and only if~$1/2 \le q \le 1$.
\begin{theorem}
    \label{th:power1n}
    Let~$q_1$ and~$q_2$ be positive constants. Suppose that~$\col = 1$ or~$\col = n$.
\begin{enumerate}
    \item\label{it:sigma1n} The function~$\res_1(X): = \|X_{-}\|_\fro^{q_1} +\|\sigma(X) -\ones\|_2^{q_2}$
          defines a local error bound for~$\stf{n,\col}_+$\!~relative
        to~$\RR^{n\times \col}$ if and only if~$q_1\le 1$ and~$q_2\le 1$, and it defines a
        global error bound if and only if~$q_1\le q_2 = 1$.
    \item\label{it:fro1n} The function~$\res_2(X): = \|X_{-}\|_\fro^{q_1} +\|X^\trs X -I_\col\|_\fro^{q_2}$
        defines a local error bound for~$\stf{n,\col}_+$\!~relative
        to~$\RR^{n\times \col}$ if and only if~$q_1\le 1$ and~$q_2\le 1$, and it defines a global error bound if and
        only if~$q_1 \le 1$ and~$1/2 \le q_2 \le 1$.
\end{enumerate}
\end{theorem}
\begin{proof}
    We consider only the case with~$\col = n$. The other case is similar.

    \ref{it:sigma1n}~Based on~\eqref{eq:lerrbdnr} and~\eqref{eq:lowerdist}, it is easy to check
    that~$\res_1$ defines a local error bound for~$\stf{n,n}_+$ relative to~$\RR^{n\times n}$
    if and only if~$q_1 \le 1$ and~$q_2\le 1$.
    Hence we only need to consider the global error bound.

    Suppose that~$q_1 \le q_2 = 1$. Let us show that
    \begin{equation}
        \label{eq:gerrbdqn}
        \dist(X,\, \stf{n,n}_+)
        \;\le\; 9n\left( \|X_{-}\|_\fro^{q_1} + \|\sigma(X) -\ones\|_2 \right)
        \;=\; 9n\res_1(X)
    \end{equation}
    for~$X\in\RR^{n\times n}$.
    If~$\|X_{-}\|_\fro \le 1$, then~\eqref{eq:gerrbdqn} follows from~\eqref{eq:gerrbdnr}.
    When~$\|X_{-}\|_\fro > 1$,
    \begin{equation}
        \nonumber
        \dist(X, \stf{n, n}_+)
        \;\le\; \dist(X, \stf{n, n})  + 2\sqrt{n}
        \;\le\; \|\sigma(X) -\ones\|_2 + 2\sqrt{n}\|X_{-}\|_\fro^{q_1},
    \end{equation}
    which validates~\eqref{eq:gerrbdqn} again.
    Hence~$\res_1$ defines a global error bound for~$\stf{n,n}_+$ relative to~$\RR^{n\times n}$.

    Now suppose that~$\res_1$ defines a global error bound for~$\stf{n,n}_+$ relative
    to~$\RR^{n\times n}$. Then it also defines a local error bound, implying~$q_1 \le 1$ and~$q_2\le 1$.
    Consider a sequence~$\{X_k\}\subset \RR^{n\times n}_+$ such that~$X_k^\trs X_k = kI_n$
    for each~$k\ge 1$. Then
    \begin{align*}
                \dist(X_k,\, \stf{n,n}_+) \;\ge\;
        \dist(X_k, \, \stf{n, n}) \;=\; \|\sigma(X_k) -\ones\|_2 \;=\; [\res_1(X_k)]
        ^{\frac{1}{q_2}} \;\to\; \infty.
    \end{align*}
    By assumption,~$\dist(X_k,\, \stf{n,n}_+)\le \nu \res_1(X_k)$ for each~$k\ge 1$ with a constant~$\nu$.
    Hence we know~$q_2 \ge 1$. To conclude, we have~$q_1 \le q_2 = 1$.
    The proof for~\ref{it:sigma1n} is complete.

    \ref{it:fro1n}~Based on~\eqref{eq:lerrbdnr},~\eqref{eq:lowerdist},
    and the fact that~$\|\sigma(X) -\ones\|_2\le\|X^\trs X - I_n\|_\fro$ (Lemma~\ref{lem:distm}),
    it is easy to check that~$\res_2$ defines a local error bound for~$\stf{n,n}_+$ relative
    to~$\RR^{n\times n}$ if and only if~$q_1 \le 1$ and~$q_2 \le 1$.
     Hence we consider only the global error bound.

    Suppose that~$q_1 \le 1$ and~$1/2 \le q_2 \le 1$. We will show that
    \begin{equation}
        \label{eq:res2bd}
        \dist(X, \stf{n, n}_+)
        \;\le\; 9n\left( \|X_{-}\|_\fro^{q_1} + \|X^\trs X -I_n\|_\fro^{q_2} \right)
        \;=\; 9n\res_2(X)
    \end{equation}
    for~$X\in\RR^{n\times n}$.
    If~$\|X^\trs X - I_n\|_\fro \le 1$, then~\eqref{eq:res2bd} holds
    because of~\eqref{eq:gerrbdqn} and the fact that $\|\sigma(X) -\ones\|_2 \le \|X^\trs X - I_n\|_\fro$.
    When~$\|X^\trs X - I_n\|_\fro > 1$,
    \begin{equation}
        \nonumber
        \begin{split}
        \dist(X, \stf{n, n}_+) \;\le\; &\dist(X, \stf{n, n})  + 2\sqrt{n}  \\
        \;\le\;&  n^{\frac{1}{4}} \|X^\trs X - I_n\|_\fro^{\frac{1}{2}} + 2\sqrt{n}  \|X^\trs X - I_n\|_\fro^{q_2} \\
        \;\le\;&  (n^{\frac{1}{4}} + 2\sqrt{n})\|X^\trs X - I_n\|_\fro^{q_2},
        \end{split}
    \end{equation}
    justifying~\eqref{eq:res2bd} again, where the second inequality applies Lemma~\ref{lem:distm}.
    Hence~$\res_2$ defines a global error bound for~$\stf{n,n}_+$ relative to~$\RR^{n\times n}$.

    Now suppose that~$\res_2$ defines a global error bound for~$\stf{n,n}_+$ relative
    to~$\RR^{n\times n}$. Then~$q_1 \le 1$ and~$q_2 \le 1$, as~$\res_2$ also defines a local
    error bound.
    Consider again a sequence~$\{X_k\}\subset \RR^{n\times n}_+$ such that~$X_k^\trs X_k = kI_n$
    for each~$k\ge 1$. Then
    \begin{gather*}
        \dist(X_k,\, \stf{n,n}_+) \ge \|\sigma(X_k) -\ones\|_2 =  (\sqrt{k}-1) \sqrt{n},\\
        \res_2(X_k) = \|X_k^\trs X_k -I_n\|_\fro^{q_2} = [(k-1)\sqrt{n}]^{q_2}.
    \end{gather*}
    By assumption,~$\dist(X_k,\, \stf{n,n}_+)\le \nu \res_2(X_k)$ for each~$k\ge 1$ with a constant~$\nu$.
  {Hence we have}~$q_2 \ge 1/2$.
    The proof for~\ref{it:fro1n} is complete.
\end{proof}

\subsection{Error bounds with~$1 \le \col \le n$}
\label{ssec:errbdr}

Now we shift our attention to the general case with~$1 \le \col \le n$.
Given previous bounds for~$\col\in\{1, n\}$, we are particularly interested in the situation where~$1 < \col < n$.

We will first prove a local error bound for~$\stf{n,\col}_+$ relative to~$\RR^{n\times \col}_+$
as detailed in Proposition~\ref{prop:lerrbdr+}.
This bound will play a role similar to what observation~\eqref{eq:keyeq1} and
Proposition~\ref{prop:lerrbdnr+} do in the cases of~$\col = 1$ and~$\col=n$, respectively.
To simplify its proof, we start with the following lemma.
\begin{lemma}
    \label{lem:xyz}
    For any matrix~$X\in\RR^{n\times \col}_+$, there exists a matrix~$Y\in\RR^{n\times \col}_+$ such
    that~$Y^\trs Y$ is diagonal and
    \begin{equation}
        \label{eq:xyz}
        \max\big\{\|x^j - y^j\|_2,\; \big|\|y^j\|_2 - 1\big|\big\}
        \;\le\; \|z^j\|_1^{\frac{1}{2}}
        \quad \text{for each} \quad j \in\{1, \dots, \col\},
    \end{equation}
    where~$x^j$,~$y^j$, and~$z^j$ denote the~$j$th column of~$X$,~$Y$, and~$Z =  X^\trs X - I_\col$,
    respectively.
\end{lemma}
\begin{proof}
    Define~$l_i$~($1\le i\le n$) and~$Y$ as in the proof of Proposition~\ref{prop:lerrbdnr+}.
    Since~$Y^\trs Y$ is diagonal as mentioned before, it suffices to establish~\eqref{eq:xyz} for this~$Y$.

    Fix an index~$j \in \{1, \dots, \col\}$.
    Recalling that~$ 0 \le X_{i, j} \le X_{i, l_i}$ for each~$i \in \{1, \dots, n\}$, we have
    \begin{equation}
        \label{eq:xdy2}
        \begin{split}
        \|x^j - y^j\|_2^2
        \;=\;
        &
        \sum_{i=1}^n X_{i,j}^2 \ind(j \neq l_i)
        \\
        \;\le\;
        &
         \sum_{i=1}^n X_{i,l_i}X_{i,j} \ind(l_i \neq j)
        \\
        \;\le\;
        &
        \sum_{l = 1}^\col \left(\sum_{i=1}^n X_{i,l}X_{i,j}\right)  \ind(l \neq j).
        \end{split}
    \end{equation}
    Since~$X^\trs X$ and~$Z$ have the same off-diagonal entries, inequality~\eqref{eq:xdy2} yields
    \begin{equation}
        \label{eq:xdy2z}
        \|x^j - y^j\|_2^2 \;\le\; \sum_{l = 1}^\col |Z_{l, j}| \ind(l \neq j)  \;=\; \|z^j\|_1 - |Z_{j,j}|.
    \end{equation}
    It remains to prove~$\big|\|y^j\|_2 - 1\big|^2 \le \|z^j\|_1$.
    To this end, note that
    \begin{equation}
        \label{eq:1my2}
        \big|\|y^j\|_2 - 1\big|^2
        \,\le\, \big|\|y^j\|_2^2 - 1\big|
        \,\le\,
        \big|\|x^j\|_2^2 - 1\big| + \|x^j - y^j\|_2^2,
    \end{equation}
    where the first inequality uses the fact that~$|t-1|^2 \le |t^2-1|$ for any~$t\ge 0$, and the
    second one is because~$\|x^j\|_2^2 - \|y^j\|_2^2 = \|x^j - y^j\|_2^2$ due to the special
    construction~\eqref{eq:y} of~$Y$.
    Since~$\|x^j\|_2^2 - 1 = Z_{j,j}$, we can combine~\eqref{eq:xdy2z} and~\eqref{eq:1my2} to obtain
    \begin{equation}
        \nonumber
        \big|\|y^j\|_2 - 1\big|^2 \;\le\;
        \big|\|x^j\|_2^2 - 1\big| +  \big(\|z^j\|_1 - |Z_{j,j}|\big) \;=\; \|z^j\|_1.
    \end{equation}
    The proof is complete.
\end{proof}

\begin{remark}
    \label{rem:jmwc}
    As mentioned earlier, the matrix~$Y$ in the proof of Lemma~\ref{lem:xyz} is the rounding
    matrix in~\cite[Procedure~3.1]{Jiang_Meng_Wen_Chen_2023}.
    Inequality~\eqref{eq:xdy2} is essentially the second inequality
    in Case II of the proof for~\cite[Lemma~3.2]{Jiang_Meng_Wen_Chen_2023}.
    The columns of~$X$ are assumed to be normalized in~\cite{Jiang_Meng_Wen_Chen_2023}, but
    such an assumption has no effect on this inequality.
\end{remark}

\begin{proposition}
    \label{prop:lerrbdr+}
    For any matrix~$X\in \RR^{n\times \col}_+$, if~$\|\sigma(X) -\ones\|_2 < 1/(3\sqrt{\col})$, then
    \begin{equation}
        \label{eq:lerrbdr+}
        \dist(X,\, \stf{n,\col}_+) \;\le\; 2\sqrt{\frac{7\col}{3}} \|\sigma(X) -\ones\|_2^{\frac{1}{2}}.
    \end{equation}
\end{proposition}
\begin{proof}
    Let~$Y$ and~$Z$ be the matrices specified in  Lemma~\ref{lem:xyz}. Then~\eqref{eq:xyz} leads to
    \begin{equation}
        \label{eq:xtoy}
        \|X-Y\|_\fro^2 \;=\; \sum_{j=1}^\col \|x^j - y^j \|_2^2  \;\le\;  \sum_{j=1}^\col \|z^j\|_1
        \;=\;\|Z\|_{\ell_1}.
    \end{equation}
    Since~$Y^\trs Y$ is diagonal, the entries of~$\sigma(Y)$ are
    \mbox{$\|y^1\|_2$,~$\dots$,~$\|y^\col\|_2$}. Thus~\eqref{eq:xyz} also provides
    \begin{equation}
    \label{eq:1msy}
    \|\sigma(Y) - \ones\|_2^2 \;=\; \sum_{j=1}^\col (\|y^j\|_2 - 1)^2 \;\le\; \sum_{j=1}^\col
    \|z^j\|_1 \;=\; \|Z\|_{\ell_1}.
    \end{equation}
    Comparing~\eqref{eq:lerrbdr+} with~\eqref{eq:xtoy}--\eqref{eq:1msy},
    we only need to prove that~$\|\sigma(X) -\ones\|_2 <\! 1/(3\sqrt{\col})$ ensures%
     \begin{equation}
         \label{eq:z1sigmax} \|Z\|_{\ell_1} \;\le\; \frac{7\col}{3} \|\sigma(X) -\ones\|_2
     \end{equation}
     and
    \begin{equation}
        \label{eq:dys} \dist(Y, \, \stf{n, \col}_+) \;=\; \|\sigma(Y) - \ones\|_2.
    \end{equation}
    Since~$\|Z\|_{\ell_1}=\sum_{i=1}^n\sum_{j=1}^\col|Z_{ij}|  \le \col\|Z\|_\fro$, inequality~\eqref{eq:z1sigmax} is a direct consequence of
    \begin{equation}
        \label{eq:zfsigmax}
        \|Z\|_\fro \;=\; \|X^\trs X -I_\col\|_\fro
        \;\le\; (\|X\|_2 + 1)\|\sigma(X) -\ones\|_2
        \;\le\; \frac{7}{3} \|\sigma(X) -\ones\|_2,
    \end{equation}
    where the last inequality is because~$\|X\|_2 \le\|\sigma(X) -\ones\|_2 + 1 < 4/3$.
    Meanwhile, inequality~\eqref{eq:zfsigmax} also leads to
    \[
        \|z^j\|_1 \;\le\; \sqrt{\col}\|Z\|_\fro \;\le\; \frac{7\sqrt{\col}}{3} \|\sigma(X) -\ones\|_2 \;<\; 1
        \quad \text{for each} \quad j\in\{1, \dots, \col\}.
    \]
    Therefore, inequality~\eqref{eq:xyz} implies that~$Y$ does not contain any zero column.
    Hence the diagonal entries of~$Y^\trs Y$ are all positive, which
    ensures the nonsingularity of this matrix since it is diagonal.
    Thus Lemma~\ref{lem:diag} yields~\eqref{eq:dys}. The proof is complete.
\end{proof}

Now we are ready to establish a local error bound for~$\stf{n,\col}_+$ relative to~$\RR^{n\times \col}$.
\begin{theorem}
    \label{th:lerrbd}
    For any matrix~$X\in \RR^{n\times \col}$, if
    $\|X_{-}\|_\fro + \|\sigma(X) - \ones\|_2 < 1/(3\sqrt{\col})$, then
    \begin{equation}
        \label{eq:lerrbd}
        \dist(X,\, \stf{n,\col}_+) \;\le\;
        4\sqrt{\col}  \left(\|X_{-}\|_\fro^{\frac{1}{2}}  + \|\sigma(X) -\ones\|_2^{\frac{1}{2}} \right).
    \end{equation}
\end{theorem}
\begin{proof}
    According to Lemma~\ref{lem:vonNeumann},
\begin{equation}
    \nonumber
    \|\sigma(X_+) -\ones\|_2 \;\le\; \|X_{-}\|_\fro + \|\sigma(X) -\ones\|_2.
\end{equation}
Thus~$\|\sigma(X_+) -\ones\|_2 < 1/(3\sqrt{\col})$ by assumption,
and hence Proposition~\ref{prop:lerrbdr+} implies
\begin{equation}
    \label{eq:x+tos}
    \begin{split}
    \dist(X_+,\, \stf{n,\col}_+)
    \;\le\; &2\sqrt{\frac{7\col}{3}}  \left(\|X_{-}\|_{\fro}^{\frac{1}{2}} + \|\sigma(X) - \ones\|_2^{\frac{1}{2}}\right).
    \end{split}
\end{equation}
On the other hand, since~$\|X_{-}\|_{\fro} < 1/(3\sqrt{\col})$, it holds that
\begin{equation}
    \label{eq:xtox+}
    \|X-X_+\|_\fro \;=\; \|X_{-}\|_\fro
    \;\le\; \frac{1}{\sqrt{3}\col^{\frac{1}{4}}} \|X_{-}\|_{\fro}^{\frac{1}{2}}
    \;\le\; \sqrt{\frac{\col}{3}} \|X_{-}\|_{\fro}^{\frac{1}{2}}.
\end{equation}
Inequality~\eqref{eq:lerrbd} follows from~\eqref{eq:x+tos} and~\eqref{eq:xtox+} because~$2\sqrt{7/3} + 1/\sqrt{3} < 4$.
\end{proof}

    Theorem~\ref{th:lerrbd} presents only a local error bound.
    Indeed,~$\|X_{-}\|_{\fro}^{\frac{1}{2}} + \|\sigma(X)-\ones\|_2^{\frac{1}{2}}$
    does not define a global error
    bound for~$\stf{n,\col}_+$ relative to~$\RR^{n\times \col}$,
    which will be explained later by Theorem~\ref{th:power}.
    To have a global error bound, we need to replace the term~$\|\sigma(X)-\ones\|_2$
    with~$\|X^\trs X - I_\col\|_\fro$ as in the following theorem. 
    \begin{theorem}\label{th:errbd}~For any matrix~$X\in \RR^{n\times \col}$, we have
    \begin{equation}
        \label{eq:gerrbd}
        \dist(X,\, \stf{n,\col}_+) \;\le\;
        5\col^{\frac{3}{4}}  \left(\|X_{-}\|_\fro^{\frac{1}{2}} + \|X^\trs X - I_\col\|_\fro^{\frac{1}{2}} \right).
    \end{equation}
    Moreover, if~$\|X_{-}\|_\fro + \|X^\trs X -I_\col\|_\fro < 1/(3\sqrt{\col})$, then
    \begin{equation}
        \label{eq:lerrbdf}
        \dist(X,\, \stf{n,\col}_+) \;\le\;
        4\sqrt{\col}  \left(\|X_{-}\|_\fro^{\frac{1}{2}}  + \|X^\trs X - I_\col\|_\fro^{\frac{1}{2}} \right).
    \end{equation}
\end{theorem}
\begin{proof}
    Recall that~$\|\sigma(X) - \ones\|_2 \le \|X^\trs X - I_\col\|_\fro$~(Lemma~\ref{lem:distm}).
    Thus~\eqref{eq:lerrbdf} is a direct consequence of Theorem~\ref{th:lerrbd}
    when~$\|X_{-}\|_\fro + \|X^\trs X -I_\col\|_\fro < 1/(3\sqrt{\col})$.

    Now we prove~\eqref{eq:gerrbd}. Let us assume that
\begin{equation}
    \nonumber
    \|X_{-}\|_{\fro} + \|X^\trs X -I_\col\|_\fro \;\ge\; \frac{1}{3\sqrt{\col}},
\end{equation}
as~\eqref{eq:gerrbd} is already justified by~\eqref{eq:lerrbdf} when this inequality does not hold.
Under this assumption,
\begin{equation}
    \label{eq:distbigsq}
        \|X_{-}\|_{\fro}^\frac{1}{2} + \|X^\trs X -I_\col\|_\fro^{\frac{1}{2}} \;\ge\;
    \frac{1}{\sqrt{3}\col^{\frac{1}{4}}}.
\end{equation}
Noting that the diameter of~$\stf{n,\col}$ is~$2\sqrt{\col}$, we then have
\begin{equation}
    \label{eq:nonlocal}
\begin{split}
    \dist(X,\, \stf{n,\col}_+) \;\le\; & \dist(X,\, \stf{n,\col}) + 2\sqrt{\col}  \\
    \;\le\; &
    \col^{\frac{1}{4}} \|X^\trs X -I_\col\|_\fro^{\frac{1}{2}}
    + 2\sqrt{3} \col^{\frac{3}{4}}
    \left( \|X_{-}\|_{\fro}^{\frac{1}{2}} + \|X^\trs X -I_\col\|_\fro^{\frac{1}{2}} \right)\\
    \;\le\; & 5\col^{\frac{3}{4}}  \left(\|X_{-}\|_\fro^{\frac{1}{2}} + \|X^\trs
    X - I_\col\|_\fro^{\frac{1}{2}} \right),
\end{split}
\end{equation}
where the second inequality applies Lemma~\ref{lem:distm} and~\eqref{eq:distbigsq}.
\end{proof}

Recently, Theorem~\ref{th:errbd} has been used in~\cite{Anthony1,Anthony2} to establish error bounds
for~$\dist(X,\, \stf{n,\col}_+)$ for~$X$ in the  unit ball of  spectral norm, i.e.,~$\{X\in
\RR^{n\times \col} : \|X\|\le 1\}$. See~(31) in~\cite{Anthony1}.

\subsection{Tightness of the error bounds when~$1 < \col < n$}
\label{ssec:tight}

The following proposition shows that the bounds presented in
Theorems~\ref{th:lerrbd} and~\ref{th:errbd} are tight up to multiplicative constants
when~$1<\col < n$,
no matter whether~$X$ belongs
to~$\stf{n,\col}$,~$\RR^{n\times \col}_+$, or neither of them.
Consequently, the error bounds~(\ref{eq:error1})--(\ref{eq:error3}) cannot hold with~$q>1/2$
when~$1< \col < n$.
\begin{proposition}
    \label{prop:tight}
    Suppose that~$1 < \col <  n$.
    \begin{enumerate}
        \item\label{it:terrbds}
            There exists a sequence~$\{X_k\}\subset\stf{n, \col}\setminus\RR^{n\times \col}_+$ such
            that~$(X_k)_{-} \to 0$ and
            \begin{equation}
                \label{eq:terrbds}
                \dist(X_k,\, \stf{n,\col}_+) \;\ge\; \frac{1}{\sqrt{2}} \|(X_k)_{-}\|_{\fro}^{\frac{1}{2}}.
            \end{equation}
        \item\label{it:terrbdr+}
            There exists a sequence~$\{X_k\}\subset\RR^{n\times \col}_+\setminus\stf{n,\col}$ such
            that~$X_k^\trs X_k \to I_\col$ and
            \begin{equation}
                \label{eq:terrbdr+}
            \dist(X_k,\, \stf{n,\col}_+) \;\ge\; \frac{1}{\sqrt{2}}\| X_k^\trs X_k - I_\col \|_{\fro}^{\frac{1}{2}}.
            \end{equation}
        \item\label{it:terrbdr}
            There exists a sequence~$\{X_k\}\subset\RR^{n\times \col} \setminus(\RR^{n\times \col}_+\cup \stf{n,\col})$
            such that~$(X_k)_{-}\to 0$, $X_k^\trs X_k \to I_\col$,  and
            \begin{equation}
                \label{eq:terrbdr}
            \dist(X_k,\, \stf{n,\col}_+)
            \;\ge\; \frac{1}{\sqrt{2}+1}\left(\|({X_k})_{-}\|_\fro^{\frac{1}{2}} + \| X_k^\trs X_k - I_\col \|_{\fro}^{\frac{1}{2}} \right).
            \end{equation}
    \end{enumerate}
\end{proposition}
\begin{proof}
    Take a sequence~$\{\varepsilon_k\}\subset (0,1/2)$ that converges to~$0$.
    For each~$k \ge 1$, let~$X_k\in\RR^{n\times \col}$ be a matrix such that its first 3 rows are
    \vspace{9pt}
  \begin{equation}
        \nonumber
 \begin{bmatrix}
           \varepsilon_k & \varepsilon_k& \overmat{\col-2}{0 & \dots  & 0}\\
            a_k & b_k & 0 & \dots  & 0 \\
            c_k & d_k & 0 & \dots  & 0
        \end{bmatrix}
    \end{equation}
    with~$a_k$,~$b_k$,~$c_k$,~$d_k$ being specified later,
    its 4th to~$(\col+1)$th rows are the last~$\col-2$ rows of~$I_\col$~(if~$\col \ge 3$), and its other
    rows are zero~(if any).
    In addition, let~$\bar{X}_k$ be a projection of~$X_k$ onto~$\stf{n,\col}_+$. Then the first
    row of~$\bar{X}_k$ contains at most one nonzero entry according to Lemma~\ref{lem:row}. Hence
    \begin{equation}
        \label{eq:distx}
        \dist(X_k,\,\stf{n,\col}_+) \;=\; \|X_k - \bar{X}_k\|_\fro \;\ge\; \varepsilon_k.
    \end{equation}
    Moreover, it is clear that~$(X_k)_{-}\to 0$ and~$X_k^\trs X_k \to I_\col$ if
    \begin{equation}
        \label{eq:abcd}
        a_k \to 1, \quad b_k \to 0, \quad c_k \to 0, \quad \text{and} \quad d_k\to 1.
    \end{equation}
    In the sequel, we will configure~$a_k$,~$b_k$,~$c_k$, and~$d_k$ subject to~\eqref{eq:abcd}
    so that~$\{X_k\}$ validates~\ref{it:terrbds},~\ref{it:terrbdr+}, and~\ref{it:terrbdr} one by one.

    \ref{it:terrbds}~Define
\begin{equation}
    \nonumber
    a_k \;=\; \sqrt{1-\varepsilon_k^2}, \quad b_k = -\frac{\varepsilon_k^2}{a_k}, \quad c_k = 0,
    \quad \text{and} \quad d_k = \sqrt{1-\varepsilon_k^2 - b_k^2}.
\end{equation}
Then~$X_k\in \stf{n,\col}\setminus\RR^{n\times \col}_+$.
Clearly,~$\|(X_{k})_{-}\|_\fro = {\varepsilon_k^2}/{a_k}$.
    Hence~\eqref{eq:terrbds} holds according to~\eqref{eq:distx} and the fact that~$a_k \ge
    \sqrt{1-\varepsilon_k^2} > 1/2$~(recall that~$\varepsilon_k <1/2$).

    \ref{it:terrbdr+}~Define~$a_k = d_k = 1$ and~$b_k = c_k = 0$. Then~$X_k\in \RR^{n\times \col}_+
    \setminus\stf{n,\col}$. By straightforward calculations,
\begin{equation}
    \nonumber
       \|X_k^\trs X_k - I_\col\|_{\fro} \;=\; 2\varepsilon_k^2.
\end{equation}
Thus~\eqref{eq:terrbdr+} holds according to~\eqref{eq:distx}.

\ref{it:terrbdr}~Define~$a_k \!= d_k \!= 1$,~$b_k = \!-\varepsilon_k^2$, and~$c_k = 0$.
Then~$X_k \!\in \RR^{n\times \col}\setminus(\RR^{n\times \col}_+ \cup \stf{n,\col})$.
In addition, we can calculate that
\begin{equation}
    \nonumber
    \|X_k^\trs X_k - I_\col\|_\fro
    \;=\; \sqrt{\varepsilon_k^4 + (\varepsilon_k^2 + \varepsilon_k^4)^2}  \;\le\;
    \sqrt{\varepsilon_k^4 + \left(\varepsilon_k^2+\frac{\varepsilon_k^2}{4}\right)^2}  \;\le\;   2\varepsilon_k^2
\end{equation}
and~$\|(X_k)_{-}\|_\fro = \varepsilon_k^2$. Therefore,~\eqref{eq:terrbdr} holds according to~\eqref{eq:distx}.
\end{proof}

Theorem~\ref{th:power} extends Theorems~\ref{th:lerrbd} and~\ref{th:errbd},
allowing~$\|X_{-}\|_\fro$ and~$\|\sigma(X) - \ones\|_2$ or
$\|X^\trs X - I_\col\|_\fro$ to have different exponents in the error bounds.
It specifies the precise range of these exponents in local and global error bounds
for~$\stf{n,\col}_+$ relative to~$\RR^{n\times\col}$ when~$1< \col < n$.
As we will see from~\ref{it:fro} of this and its proof,
when~$1< \col < n$, the error bound~\eqref{eq:error1} can hold if and only if~$q \le 1/2$,
whereas~\eqref{eq:error2} and~\eqref{eq:error3} can hold if and only if~$q = 1/2$.

\begin{theorem}
    \label{th:power}
    Let~$q_1$ and~$q_2$ be positive constants. Suppose that~$1 < \col < n$.
\begin{enumerate}
    \item\label{it:sigma} The function~$\res_1(X): = \|X_{-}\|_\fro^{q_1} +\|\sigma(X) -\ones\|_2^{q_2}$
          defines a local error bound for~$\stf{n,\col}_+$ relative
        to~$\RR^{n\times \col}$ if and only if~$q_1\le 1/2$ and~$q_2\le 1/2$, but it cannot define
        a global error bound no matter what values~$q_1$ and~$q_2$ take.
    \item\label{it:fro} The function~$\res_2(X): = \|X_{-}\|_\fro^{q_1} +\|X^\trs X -I_\col\|_\fro^{q_2}$
        defines a local error bound for~$\stf{n,\col}_+$ relative
        to~$\RR^{n\times \col}$ if and only if~$q_1\le 1/2$ and~$q_2\le 1/2$, and it defines a global error bound if and
        only if~$q_1 \le q_2 = 1/2$.
\end{enumerate}
\end{theorem}
\begin{proof}
   \ref{it:sigma}~Based on~\eqref{eq:lerrbd}, it is easy to check that~$\res_1$ defines a local error bound
    for~$\stf{n,\col}_+$ relative to~$\RR^{n\times \col}$ if~$q_1 \le 1/2$ and~$q_2\le 1/2$.
    Conversely, if~$\res_1$ defines a local error bound for~$\stf{n,\col}_+$ relative
    to~$\RR^{n\times \col}$, then~$q_1 \le 1/2$ and~$q_2\le 1/2$ according to~\ref{it:terrbds}
    and~\ref{it:terrbdr+} of Proposition~\ref{prop:tight}, respectively.

    Now we prove that~$\res_1$ cannot define a global error bound.
    According to what has been shown above, we assume that~$q_2 \le 1/2$, as a global error bound must be a local one.
    Consider a sequence~$\{X_k\}\subset \RR^{n\times \col}_+$ with~$\|X_k\|_\fro \to \infty$.
    Then~$\res_2(X_k) = \|\sigma(X_k)-\ones\|_2^{q_2}$, and hence
    \begin{equation*}
        \frac{\dist(X_k,\, \stf{n,\col}_+)}{\res_1(X_k)}
        \;\ge\;\frac{\|\sigma(X_k) -\ones\|_2}{\|\sigma(X_k) -\ones\|_2^{q_2}}
        \;\to\; \infty.
    \end{equation*}
   Thus~$\res_1$ cannot define a global error bound for~$\stf{n,\col}_+$ relative to~$\RR^{n\times \col}$.

    \ref{it:fro}~Similar to~\ref{it:sigma}, we can show that~$\res_2$ defines a local error bound
    for~$\stf{n,\col}_+$ relative to~$\RR^{n\times \col}$ if and only if~$q_1 \le 1/2$ and~$q_2 \le 1/2$.
    Hence we only need to consider the global error bound.

    Suppose that~$q_1 \le q_2 = 1/2$. Let us show that
    \begin{equation}
        \label{eq:gerrbdq}
        \dist(X,\, \stf{n,\col}_+)
        \;\le\; 5\col^{\frac{3}{4}}\left( \|X_{-}\|_\fro^{q_1} + \|X^\trs X -I_\col\|_\fro^{\frac{1}{2}} \right)
        \;=\; 5\col^{\frac{3}{4}}\res_2(X)
    \end{equation}
    for all~$X\in\RR^{n\times \col}$.
    If~$\|X_{-}\|_\fro \le 1$, then~\eqref{eq:gerrbdq} follows from~\eqref{eq:gerrbd}. When~$\|X_{-}\|_\fro > 1$,
    \begin{equation}
        \nonumber
        \dist(X,\,\stf{n, \col}_+)
        \;\le\; \dist(X,\,\stf{n, \col}) + 2\sqrt{\col}
         \;\le\; \col^{\frac{1}{4}}\|X^\trs X -I_\col\|_\fro^{\frac{1}{2}} + 2\sqrt{\col} \|X_{-}\|_\fro^{q_1}
         \;\le\; 5\col^{\frac{3}{4}}\res_2(X),
    \end{equation}
    where the second inequality applies Lemma~\ref{lem:distm}.
    Hence~$\res_2$ defines a global error bound for~$\stf{n,\col}_+$ relative to~$\RR^{n\times \col}$.

    Now suppose that~$\res_2$ defines a global error bound for~$\stf{n,\col}_+$ relative
    to~$\RR^{n\times \col}$. Then it defines a local error bound, implying~$q_1 \le 1/2$
    and~$q_2 \le 1/2$.
    Similar to the proof for~\ref{it:fro1n} of Theorem~\ref{th:power1n}, by
    considering a sequence~$\{X_k\}\subset \RR^{n\times \col}_+$ such that~$X_k^\trs X_k = kI_\col$
    for each~$k\ge 1$, we can prove~$q_2 \ge 1/2$.
    The proof is complete.
\end{proof}

\subsection{Linear regularity of~$\RR^{n\times \col}_+$ and~$\stf{n,\col}$}
\label{ssec:linreg}

Before ending this section, we briefly mention that our analysis enables us to characterize the
linear regularity of~$\RR^{n\times \col}_+$ and~$\stf{n,\col}$ for~$\col\in\{1, \dots, n\}$.

A pair of sets~$\AA_1$ and~$\AA_2$ in~$\RR^{n\times \col}$ with~$\AA_1\cap\AA_2 \neq
\emptyset$ are said to be boundedly linearly regular if for any bounded
set~$\TT \subset \RR^{n\times \col}$ there exists a constant~$\gamma$ such that%
\begin{equation}
\label{eq:linreg}
\dist(X,\, \AA_1\cap \AA_2) \;\le\; \gamma \max\left\{ \dist(X,\, \AA_1),\; \dist(X,\, \AA_2) \right\}
\end{equation}
for all~$X \in \TT$, and they are linearly regular if~\eqref{eq:linreg} holds for
all~$X\in\RR^{n\times \col}$.
Linear regularity is a fundamental concept in
optimization and is closely related to error bounds. See~\cite{Ng_Yang_2004}
and~\cite[Section~8.5]{Cui_Pang_2021} for more details.
Note that we can replace the maximum in~\eqref{eq:linreg} with a summation without essentially
changing the definition of~(boundedly) linear regularity.

Proposition~\ref{prop:linreg} clarifies whether~$\RR^{n\times \col}_+$ and~$\stf{n,\col}$ are linearly regular.
\begin{proposition}\label{prop:linreg} The two sets
    $\RR^{n\times \col}_+$~and~$\stf{n,\col}$~are linearly regular if and only
    if~$\col = 1$ or~$\col = n$.%
\end{proposition}
\begin{proof}
Recall that~$\dist(X,\,\RR^{n\times\col}_+) = \|X_{-}\|_\fro$
and~$\dist(X,\,\stf{n,\col}) = \|\sigma(X) -\ones\|_2$ for~$X\in\RR^{n\times \col}$.
The ``if'' part of this proposition holds because of the global error bounds in
Theorems~\ref{th:gerrbdr1} and~\ref{th:errbdnr}. The ``only if'' part holds
because~$\|X_{-}\|_\fro + \|\sigma(X) -\ones\|_2$ does not define a global error bound
for~$\stf{n,\col}_+$ relative to~$\RR^{n\times\col}$ when~$1<\col <n$,
as we can see from~\ref{it:sigma} of Theorem~\ref{th:power}.
\end{proof}

Proposition~\ref{prop:linreg} remains true if we change ``linearly regular'' to ``boundedly linearly
regular''. The ``if'' part is weakened after this change, and the other part holds
because~$\|X_{-}\|_\fro + \|\sigma(X) -\ones\|_2$ does not define a local error bound
for~$\stf{n,\col}_+$ relative to~$\RR^{n\times\col}$ when~$1<\col <n$
according to~\ref{it:sigma} of Theorem~\ref{th:power}.

\section{Error bounds for the sign-constrained Stiefel manifold}
\label{sec:serrbd}

This section will establish the error bounds for~$\stf{n,\col}_\tsign$ based on those already
proved for~$\stf{n,\col}_+$.

\subsection{A special case} 

First, we consider the special case with
\begin{equation*}
    \posset = \{1, \dots, \col_1\} ~\text{ and }~ \negset = \emptyset,
\end{equation*}
where~$\col_1\in \{1, \dots, \col\}$.
Define~$\col_2 = \col-\col_1$ henceforth.
In this case,~$\stf{n,\col}_\tsign$ reduces to
\begin{equation}
    \label{eq:stfr1}
    \stf{n,\col}_{\col_1,+} \;\coloneqq\;
    \left\{ X = (X_1, X_2) \mathrel{|} X_1\in\RR^{n\times \col_1}_+, ~X_2\in\RR^{n\times \col_2},
    ~X^\trs X = I_\col \right\},
\end{equation}
with~$\stf{n,\col}_{\col_1,+}$ being~$\stf{n,\col}_+$ if~$\col_1 = \col$.

Note that the results established in Sections~\ref{sec:lemmas} and~\ref{sec:errbd}
are still valid when~$\col$ is replaced with~$\col_1$ or~$\col_2$. In the sequel, we will
apply these results directly without restating this fact.

\begin{lemma}
\label{lem:sigerrbd}
Suppose that~$\col_1<\col$.
Consider matrices~$Y_1 \in \RR^{n\times\col_1}$ and~$Y_2\in \RR^{n\times \col_2}$. If~$Y_1^\trs Y_2 = 0$,
then there exists a matrix~$Z$ that is a projection of~$Y_2$ onto~$\stf{n,\col_2}$ and
satisfies~$Y_1^\trs Z= 0$.
\end{lemma}
\begin{proof}
Define~$k = n - \rank(Y_1)$. Take a matrix~$V\in\stf{n, k}$ such that~$\range(V)$ is the orthogonal
complement of~$\range(Y_1)$ in~$\RR^n$. Since~$k \ge \col -\col_1 = \col_2$,
the matrix~$V^\trs Y_2\in\RR^{k\times \col_2}$ has a polar decomposition~$UP$
with~$U\in\stf{k, \col_2}$ and~$P\in\RR^{\col_2\times \col_2}$, the latter being positive semidefinite.
Define~$Z= VU \in\RR^{n\times \col_2}$.
Then
\[
    Z P = VUP =VV^\trs Y_2 = Y_2,
\]
where the last equality holds because~$\range(Y_2) \subset \range(V)$ according to~$Y_1^\trs Y_2 = 0$,
and~$VV^\trs$ is the orthogonal projection onto~$\range(V)$.
Besides,~$Z^\trs Z = U^\trs V^\trs VU = I_{\col_2}$. Thus~$Z P$ is a polar decomposition of~$Y_2$.
Hence~$Z$ is a projection of~$Y_2$ onto~$\stf{n,\col_2}$ by Lemma~\ref{lem:FanHoffman}.
Moreover,~$Y_1^\trs Z = Y_1^\trs VU = 0$.
\end{proof}

Note that~$\stf{n,\col}_{\col_1,+}$ can also be formulated as
\[
    \stf{n,\col}_{\col_1,+}
    \;=\; \left\{ (X_1, X_2) \mathrel{|} X_1\in\stf{n, \col_1}_+, ~X_2\in\stf{n, \col_2},
    ~X_1^\trs X_2 = 0\right\}.
\]
This formulation motivates us to develop the following lemma, which provides a global error bound
for~$\stf{n,\col}_{\col_1,+}$ relative to~$\RR^{n\times \col}$.
\begin{lemma}\label{lem:X2signStf}
    Suppose that~$\col_1 < \col$.
For any matrix~$X = (X_1, X_2)$ with~$X_1\in \RR^{n\times \col_1}$ and~$X_2\in \RR^{n\times \col_2}$,
we~have%
\begin{equation}\label{X2signStf}
    \dist(X,\, \stf{n,\col}_{\col_1,+})\;\le\;
    (2\|X_2\|_2+1)\dist(X_1,\, \stf{n,\col_1}_+) + \dist(X_2,\, \stf{n,\col_2}) + 2\|X_1^\trs X_2\|_\fro.
\end{equation}
\end{lemma}

\begin{proof}
Let~$Y_1$ be a projection of~$X_1$ onto~$\stf{n, \col_1}_+$
and~$Y_2 = (I_{n} - Y_1Y_1^\trs)X_2\in\RR^{n\times \col_2}$.
Then~$Y_1^\trs Y_2 = 0$.
By Lemma~\ref{lem:sigerrbd}, there exists a matrix~$Z$ that is a projection of~$Y_2$
onto~$\stf{n, \col_2}$ with~$Y_1^\trs Z = 0$.
Define~$\bar{X} = (Y_1, Z)$, which lies in~$\stf{n, \col}_{\col_1,+}$.
Let us estimate~$\|X - \bar{X}\|_\fro$. It is clear that
\begin{equation}
    \nonumber
\begin{split}
\|X - \bar{X}\|_\fro
\le\; & \|(X_1, X_2) - (Y_1, Y_2)\|_\fro + \| (Y_1, Y_2) - (Y_1, Z)\|_\fro \\
\;\le\; & \|X_1 - Y_1\|_\fro + \|X_2 - Y_2\|_\fro + \|Y_2 - Z\|_\fro.
\end{split}
\end{equation}
Since~$\|Y_2-Z\|_\fro = \|\sigma(Y_2) - \ones\|_2$ (Lemma~\ref{lem:distm}) and
$\|\sigma(X_2)-\sigma(Y_2)\|_2\le\|X_2-Y_2\|_\fro$ (Lemma~\ref{lem:vonNeumann}),
it holds that~$\|Y_2-Z\|\le\|\sigma(X_2) - \ones\|_2 + \|X_2 - Y_2 \|_\fro$.
Therefore,
\begin{equation}
    \label{eq:XmZ}
    \|X-\bar{X}\|_\fro \;\le\; \|X_1-Y_1\|_\fro + \|\sigma(X_2) - \ones\|_2 + 2\|X_2 - Y_2\|_\fro.
\end{equation}
Meanwhile, recalling that~$Y_2=(I_n - Y_1Y_1^\trs) X_2$ and~$Y_1\in\stf{n,\col_1}$,
we have
\begin{equation}
    \label{eq:X2mY}
    \begin{split}
    \|X_2 - Y_2\|_\fro =\; & \|Y_1Y_1^\trs X_2\|_\fro \;= \;\|Y_1^\trs X_2\|_\fro
    \;\le\;\|(Y_1-X_1)^\trs X_2\|_2 + \|X_1^\trs X_2\|_\fro.
    \end{split}
\end{equation}
Plugging~\eqref{eq:X2mY} into~\eqref{eq:XmZ}
while noting~$\|(Y_1-X_1)^\trs X_2\|\le\|X_1-Y_1\|_\fro\|X_2\|_2$, we obtain%
\begin{equation*}
    \|X-\bar{X}\|_\fro \;\le\; (2\|X_2\|_2 + 1)\|X_1-Y_1\|_\fro + \|\sigma(X_2) - \ones\|_2 + 2\|X_1^\trs X_2\|_\fro.
\end{equation*}
This implies inequality~\eqref{X2signStf}, because~$\|X_1-Y_1\|_\fro =\dist(X_1,\,\stf{n,\col_1}_+)$
by the definition of~$Y_1$
and~$\|\sigma(X_2) - \ones\|_2 =\dist(X_2,\,\stf{n,\col_2})$ by Lemma~\ref{lem:distm}.
\end{proof}

In light of Lemma~\ref{lem:X2signStf}, we can establish error bounds for~$\stf{n,\col}_{\col_1,+}$
using those for~$\stf{n,\col}_{+}$, as will be done in Propositions~\ref{prop:errbds1}
and~\ref{prop:errbds}. To this end, it is useful to note for any matrix~$X=(X_1, X_2)$ that
\begin{equation}
    \label{eq:sumdist}
    \|X^\trs X -I_{\col}\|_\fro \;\ge\;
    \max\left\{\|X_1^\trs X_1-I_{\col_1}\|_\fro,\;
    \|X_2^\trs X_2-I_{\col_2}\|_\fro,\;
\sqrt{2}\|X_1^\trs X_2\|_\fro\right\}.
\end{equation}

\begin{proposition}
    \label{prop:errbds1}
    For any matrix~$X \in \RR^{n\times \col}$ with~$x_1$ being its first column, we~have
    \begin{equation}
        \label{eq:gerrbds1}
        \dist(X,\, \stf{n,\col}_{1,+}) \;\le\;
        7\sqrt{\col}\left(\|(x_1)_{-}\|_2 + \|X^\trs  X - I_\col\|_\fro\right).
    \end{equation}
        Moreover, if~$\|X^\trs X -I_\col\|_\fro < 1/3$, then
    \begin{equation}
        \label{eq:lerrbds1}
        \dist(X,\, \stf{n,\col}_{1,+}) \;\le\;
        7\left( \|(x_1)_{-}\|_2 + \|X^\trs  X - I_\col\|_\fro \right).
    \end{equation}
\end{proposition}
\begin{proof}
    If~$\col = 1$, then~\eqref{eq:gerrbds1} and~\eqref{eq:lerrbds1} hold because of
    Theorem~\ref{th:gerrbdr1}. Hence we suppose that~$\col > 1$ in the sequel.
     We first assume~$\|X^\trs  X - I_\col\|_\fro < 1/3$ and establish~\eqref{eq:lerrbds1}.
    Let~$X_2$ be the matrix containing the last~$\col-1$ columns of~$X$.
 According to Theorem~\ref{th:gerrbdr1} and Lemma~\ref{lem:distm},
    \begin{gather}
        \label{eq:x1} \dist(x_1,\,\stf{n, 1}_+)  
        \;\le\; 2\|(x_1)_-\|_2 +  \big|x_1^\trs x_1 - 1\big|,\\
        \label{eq:x2} \dist(X_2,\,\stf{n, \col -1})  \;\le\; \|X_2^\trs X_2-I_{\col-1}\|_\fro.
    \end{gather}
    Plugging~\eqref{eq:x1} and~\eqref{eq:x2} into Lemma~\ref{lem:X2signStf} while
    noting~\eqref{eq:sumdist}, we have
\begin{equation}
    \nonumber
    \begin{split}
    \dist(X,\, \stf{n,\col}_{1, +})
    & \;\le\; (2\|X_2\|_2+1)\cdot 2\|(x_1)_-\|_2 + \big[(2\|X_2\|_2+1)+1+\sqrt{2}\big]\|X^\trs X-I_{\col}\|_\fro\\
    & \;\le\; 7\left(\|(x_1)_-\|_2 + \|X^\trs X-I_{\col}\|_\fro\right),
    \end{split}
\end{equation}
where the second inequality uses the fact that~$\|X_2\|_2 \le \|X^\trs X \|_\fro^{\frac{1}{2}} \le 2/\sqrt{3}$.

To prove~\eqref{eq:gerrbds1}, we now only need to focus on the case
with~$\|X^\trs X - I_\col\|_\fro \ge 1/3$. In this case,
\begin{equation}
    \nonumber
    \dist(X, \stf{n, \col}_{1,+})
    \;\le\; \dist(X, \stf{n,\col}) + 2\sqrt{\col}
    \;\le\; \|X^\trs X - I_\col\|_\fro + 6\sqrt{\col}\|X^\trs X - I_\col\|_\fro,
\end{equation}
which implies~\eqref{eq:gerrbds1}. The proof is complete.
\end{proof}

\begin{proposition}\label{prop:errbds}
    For any matrix~$X \in\RR^{n\times \col}$ with~$X_1$ being its submatrix
    containing the first~$\col_1$ columns, we~have
    \begin{equation}
        \label{eq:gerrbds}
        \dist(X,\, \stf{n,\col}_{\col_1,+}) \;\le\;
        15\col^{\frac{3}{4}} \left(\|(X_1)_{-}\|_\fro^{\frac{1}{2}} + \|X^\trs  X - I_\col\|_\fro^{\frac{1}{2}} \right).
    \end{equation}
        Moreover, if~$\|(X_1)_{-}\|_\fro + \|X^\trs X -I_\col\|_\fro < 1/(3\sqrt{\col})$, then
    \begin{equation}
        \label{eq:lerrbds}
        \dist(X,\, \stf{n,\col}_{\col_1,+}) \;\le\;
        15\sqrt{\col} \left( \|(X_1)_{-}\|_\fro^{\frac{1}{2}} + \|X^\trs  X - I_\col\|_\fro^{\frac{1}{2}} \right).
    \end{equation}
\end{proposition}
 \begin{proof}
     If~$\col_1 = \col$, then~\eqref{eq:gerrbds} and~\eqref{eq:lerrbds} hold because of
     Theorem~\ref{th:errbd}. Hence we suppose that~$\col_1 < \col$ in the sequel.
     We first assume~$\|(X_1)_{-}\|_\fro + \|X^\trs  X - I_\col\|_\fro < 1/(3\sqrt{\col})$ and establish~\eqref{eq:lerrbds}.
     Let~$X_2$ be the matrix containing the last~$\col_2 = \col - \col_1$ columns of~$X$.
     According to~\eqref{eq:sumdist},
 our assumption implies
 \begin{gather}
     \nonumber
 \|(X_1)_{-}\|_\fro + \|X_1^\trs  X_1 - I_{\col_1}\|_\fro \;<\; \frac{1}{3\sqrt{\col_1}},
 \quad \|X_2^\trs X_2 - I_{\col_2}\|_\fro \;\le\; \frac{1}{3}.
 \end{gather}
Hence Theorem~\ref{th:errbd} and Lemma~\ref{lem:distm} yield
\begin{align}
        \label{eq:d1}
\dist(X_1,\, \stf{n,\col_1}_+) &
\;\le\; 4\sqrt{\col_1} \left(\|(X_1)_{-}\|_\fro^{\frac{1}{2}} + \|X^\trs_{1} X_{1} - I_{\col_1} \|_\fro^{\frac{1}{2}}\right),\\
        \label{eq:d2}
\dist(X_2,\, \stf{n,\col_2}) &
\;\le\;  \|X_2^\trs X_2-I_{\col_2}\|_\fro \;\le\; \frac{1}{\sqrt{3}}\|X_2^\trs X_2-I_{\col_2}\|_\fro^{\frac{1}{2}}.
\end{align}
In addition, inequality~\eqref{eq:sumdist} and our assumption also provide
\begin{equation}
  \label{eq:d3}
\|X_1^\trs X_2\|_\fro
\;\le\;  \frac{1}{\sqrt{2}}\|X^\trs X-I_{\col}\|_\fro \;\le\; \frac{1}{\sqrt{6}}\|X^\trs X-I_{\col}\|_\fro^{\frac{1}{2}}.
\end{equation}
Plugging~\eqref{eq:d1}--\eqref{eq:d3} into Lemma~\ref{lem:X2signStf}
while noting~\eqref{eq:sumdist}, we obtain
\begin{equation*}
\begin{split}
\dist(X,\, \stf{n,\col}_{\col_1,+}) \;\le\; &
\left[4\sqrt{\col_1}(2\|X_2\|_2 + 1) + \frac{1}{\sqrt{3}} + \frac{2}{\sqrt{6}} \right]
\left(\|(X_1)_{-}\|_\fro^{\frac{1}{2}} + \|X^\trs X -I_\col\|_\fro^{\frac{1}{2}}\right)\\
\le\; & 15 \sqrt{\col} \left(\|(X_1)_{-}\|_\fro^{\frac{1}{2}} + \|X^\trs X -I_\col\|_\fro^{\frac{1}{2}}
\right),
\end{split}
\end{equation*}
where the second inequality uses the fact that~$\|X_2\|_2\le\|X^\trs  X\|_\fro^{\frac{1}{2}} \le 2/\sqrt{3}$.

Now we prove~\eqref{eq:gerrbds}. By the same technique as the proof of~\eqref{eq:nonlocal}, we have
\begin{equation*}
\dist(X,\, \stf{n,\col}_{\col_1,+}) \;\le\; 5 \col^{\frac{3}{4}} \left(\|(X_1)_{-}\|_\fro^{\frac{1}{2}}
+  \|X^\trs X - I_\col\|_\fro^{\frac{1}{2}}\right)
\end{equation*}
when~$\|(X_1)_{-}\|_\fro + \|X^\trs  X - I_\col\|_\fro \ge 1/(3\sqrt{\col})$.
Combining this with~\eqref{eq:lerrbds}, we conclude that~\eqref{eq:gerrbds} is valid. The proof is
complete.
\end{proof}

\subsection{The general case}

We now present the error bounds for~$\stf{n,\col}_\tsign$, detailed in
Theorems~\ref{th:errbdsgen1}--\ref{th:errbdsgen}.
Theorems~\ref{th:errbdsgen1} and~\ref{th:errbdsgenn} can be proved using Proposition~\ref{prop:errbds1}
and Theorem~\ref{th:errbdnr}, respectively. We omit the proofs because they are essentially the same as
that of Theorem~\ref{th:errbdsgen} below.
\begin{theorem}\label{th:errbdsgen1}
    Suppose that~$|\posset| + |\negset| = 1$. For any matrix~$X\in \RR^{n\times \col}$, we~have
    \begin{equation}
        \nonumber
        \dist(X,\, \stf{n,\col}_{\tsign}) \;\le\;
        7\sqrt{\col} \left(\|(\signmat \hp X)_{-}\|_\fro+ \|X^\trs  X - I_\col\|_\fro \right).
    \end{equation}
    Moreover, if~$\|X^\trs X -I_\col\|_\fro < 1/3$, then
    \begin{equation}
        \nonumber
        \dist(X,\, \stf{n,\col}_{\tsign}) \;\le\;
        7 \left( \|(\signmat\hp X)_{-}\|_\fro + \|X^\trs  X - I_\col\|_\fro \right).
    \end{equation}
\end{theorem}

\begin{theorem}
    \label{th:errbdsgenn}
    Suppose that~$|\posset| + |\negset| =n$. For any matrix~$X\!\in\RR^{n\times n}$,~we~have%
        \begin{equation}
            \nonumber
            \dist(X,\,\stf{n,n}_\tsign) \;\le\; 9n\left( \|(\signmat \hp X)_{-}\|_\fro + \|\sigma(X) -\ones\|_2 \right).
        \end{equation}
        Moreover, if~$\|(\signmat \hp X)_{-}\|_\fro + \|\sigma(X) -\ones\|_2 < 1/(4\sqrt{n})$, then
        \begin{equation}
            \nonumber
            \dist(X,\,\stf{n,n}_\tsign) \;\le\; 8\sqrt{n}\left( \|(\signmat \hp X)_{-}\|_\fro + \|\sigma(X) -\ones\|_2 \right).
        \end{equation}
\end{theorem}

\begin{theorem}\label{th:errbdsgen}
    For any matrix~$X\in \RR^{n\times \col}$, we~have
    \begin{equation}
        \label{eq:gerrbdsgen}
        \dist(X,\, \stf{n,\col}_{\tsign}) \;\le\;
        15\col^{\frac{3}{4}} \left(\|(\signmat \hp X)_{-}\|_\fro^{\frac{1}{2}} + \|X^\trs  X - I_\col\|_\fro^{\frac{1}{2}} \right).
    \end{equation}
    Moreover, if~$\|(\signmat\hp X)_{-}\|_\fro + \|X^\trs X -I_\col\|_\fro < 1/(3\sqrt{\col})$, then
    \begin{equation}
        \label{eq:lerrbdsgen}
        \dist(X,\, \stf{n,\col}_{\tsign}) \;\le\;
        15\sqrt{\col} \left( \|(\signmat\hp X)_{-}\|_\fro^{\frac{1}{2}} + \|X^\trs  X - I_\col\|_\fro ^{\frac{1}{2}} \right).
    \end{equation}
\end{theorem}

\begin{proof}
    Let~$\QQ = \{1, \dots, \col\}\setminus (\posset \cup \negset)$.
    With~$M_{\sss{\posset}}$, $M_{\sss{\negset}}$, and~$M_{\QQ}$ being
    the submatrices of~$I_\col$ containing the columns indexed by~$\posset$, $\negset$, and~$\QQ$,
    respectively, we take
    the permutation matrix
    \[
        \Pi \;=\;  (M_{\sss{\posset}}, \; M_{\sss{\negset}}, \; M_{\sss{\QQ}}) \;\in\;
        \RR^{\col\times\col}.
    \]
    In addition, we take the diagonal matrix~$D\in\RR^{\col\times\col}$ with~$D_{j,j} = -1$
    if~$j\in\negset$ and~$D_{j,j} = 1$ otherwise.
    Define~$\col_1 = |\posset| + |\negset|$. If~$\col_1 = 0$, then~\eqref{eq:gerrbdsgen}
    and~\eqref{eq:lerrbdsgen} hold because of Lemma~\ref{lem:distm}. Hence we suppose
    that~$\col_1\ge 1$ in the sequel.

    Consider any matrix~$X\in\RR^{n\times \col}$. Let~$Y = X D \Pi$, and~$\bar{Y}$ be the projection
    of~$Y$ onto~$\stf{n,\col}_{\col_1,+}$ defined in~\eqref{eq:stfr1}.
    Set~$\bar{X} = \bar{Y}\Pi^\trs D$, which lies in~$\stf{n,\col}_\tsign$. Then
    \[
        \dist(X,\, \stf{n,\col}_\tsign) \;\le\; \|X - \bar{X}\|_\fro \;=\;\|Y\Pi^\trs D - \bar{Y}
        \Pi^\trs D\|_\fro
        \;=\;\|Y - \bar{Y}\|_\fro.
    \]
    Invoking Proposition~\ref{prop:errbds}, we have
    \[
        \|Y-\bar{Y}\|_\fro \;\le\;
        15\col^{\frac{3}{4}}
        \left( \|(Y_1)_{-}\|_\fro^{\frac{1}{2}} + \|Y^\trs  Y - I_\col\|_\fro ^{\frac{1}{2}} \right),
    \]
    where~$Y_1$ is the submatrix of~$Y$ containing the first~$\col_1$ columns. It is
    straightforward to verify that~$\|(Y_1)_{-}\|_\fro = \|(\signmat\hp X)_{-}\|_\fro$
    and~$\|Y^\trs Y - I_\col\|_\fro = \|X^\trs X - I_\col\|_\fro$. Hence we obtain~\eqref{eq:gerrbdsgen}. The
    bound~\eqref{eq:lerrbdsgen} can be established in a similar way.
\end{proof}

\section{Exact penalties for optimization on the nonnegative Stiefel manifold}
\label{sec:exactpen}

 In this section, as an application of the error bounds established in this paper,
we consider exact penalties for optimization problem~(\ref{eq:optnoc}). For simplicity, we will
focus on the special case with
\[
    \stf{n,\col}_{\signmat} = \stf{n,\col}_+,
\]
applying the bounds in Section~\ref{sec:errbd}.
Essentially the same results can be established in the general case
by exploiting the bounds in Section~\ref{sec:serrbd}.
The exact penalty results only require~(local)~Lipschitz continuity of~$F$, and hence
can be applied to nonsmooth optimization,
for example,~$F$ involving a group sparse regularization term~\cite{Xiao_Liu_Yuan_2021}.

The exactness of penalty methods for problem~\eqref{eq:optnoc}
with~$\stf{n,\col}_{\signmat} = \stf{n,\col}_+$ has been studied
in~\cite{Jiang_Meng_Wen_Chen_2023,Qian_Pan_Xiao_2023}.
In~\cite{Jiang_Meng_Wen_Chen_2023}, an error bound is established for~$\stf{n,\col}_+$ relative to the set
\[
\{X\in \RR^{n\times \col}_+ \mathrel{:} (X^\trs X)_{j,j}=1,\, j = 1, \dots, \col\},
\]
and then the bound is used to analyze a penalty method.
However, the
error bound in~\cite{Jiang_Meng_Wen_Chen_2023} cannot be used to  derive the values of~$\nu$ and~$q$
in~(\ref{eq:error1})--(\ref{eq:error3}).
In~\cite{Qian_Pan_Xiao_2023}, the authors consider
the penalty problem~(\ref{penalty1}) with~$q=1$, and show
this problem has the same global minimizers as problem~\eqref{eq:optnoc}
if each global optimal solution of~\eqref{eq:optnoc}  has no zero rows.  Our exact penalty results
only need the Lipschitz continuity of the objective function~$F$ in~\eqref{eq:optnoc}.

The error bounds~\mbox{(\ref{eq:error1})--(\ref{eq:error3})} established in this paper
enable us to have the exactness of the penalized problem
\begin{equation}
    \label{eq:nocpen}
    \min \left\{ F(X) + \penpar\left(\|X_{-}\|_{\ell_{p}}^{q_1} + \|X^\trs X - I_\col\|_{\ell_{p}}^{q_2} \right) \mathrel{:} X \in \sset\right\}
\end{equation}
for solving~\eqref{eq:optnoc}
with~$\stf{n,\col}_{\signmat} = \stf{n,\col}_+$
only under the~(local) Lipschitz continuity of function~$F$.
Here the set~$\sset \subset \RR^{n\times \col}$ is a set that contains~$\stf{n,\col}_+$,
while the parameters~$\penpar$,~$p$,~$q_1$, and~$q_2$ are all positive.
If~$p = 2$ and~$q_1 = q_2 = q$, then the penalized problem~\eqref{eq:nocpen} reduces to
problems~\eqref{penalty1} and~\eqref{penalty2}
when~$\sset$ equals~$\stf{n,\col}$ and~$\RR^{n\times \col}_+$, respectively.

Due to the equivalence between norms, it is indeed possible to establish the exactness
of~\eqref{eq:nocpen} when the entry-wise~$\ell_p$-norm is changed to other ones.
We choose to use the entry-wise~$\ell_p$-norm in~\eqref{eq:nocpen} because it is easy to evaluate.

\subsection{Exactness for Lipschitz continuous objective functions}\label{ssec:exactposorth}

Theorem~\ref{th:exactposorth} presents the exactness of problem~(\ref{eq:nocpen}) regarding global optimizers
when the objective function~$F\mathrel{:}\sset \to \RR$ is an~$L$-Lipschitz continuous function, namely
\begin{equation}
    \label{eq:lip}
    |F(X) - F(Y)| \;\le\; L \|X-Y\|_\fro
\end{equation}
for all~$X$ and~$Y$ in~$\sset$,
where~$L\in(0,\infty)$ is a Lipschitz constant of~$F$ with respect to the Frobenius norm.
Note that the global Lipschitz continuity of the objective function~$F$~(\ref{eq:lip}) is assumed on
a set~$\sset$ containing~$ \stf{n,\col}_+$.  For example,
if~$F(X)=$trace$(X^\trs A^\trs A X)$  and~$\sset=\{X \in \RR^{n\times \col}:  \|X\|_\fro\le \gamma\}$
with~$\gamma>\sqrt{\col}$, then the global Lipschitz continuity of~$F$ holds on~$\sset$ with the
Lipschitz constant~$L=2\gamma\|A\|_2^2.$  Indeed,  our theory holds even if~$F$ is undefined out of~$\sset$.
{The proof of Theorem~\ref{th:exactposorth} is standard and we include it in
Appendix~\ref{app:proof} for completeness.}
\begin{theorem}[Exact penalty~(\ref{eq:nocpen}) with~$F$ being Lipschitz continuous]\label{th:exactposorth}
    Suppose that~$\sset\subset\RR^{n\times\col}$ is a set
    containing~$\stf{n,\col}_+$,~$F\mathrel{:} \sset\to\RR$
    is an~$L$-Lipschitz function, and~$p\ge 1$ is a constant.
If~$0< q \le 1/2$ and~$\penpar > 5L\col^{\frac{3}{4}}\max\Big\{1, \, (n\col)^{\frac{p-2}{4p}}\Big\}$, then
\begin{equation}
    \nonumber
    \begin{split}
    \Argmin\{F(X) \mathrel{:} X \!\in \stf{n,\col}_+\}
    = \Argmin\Big\{F(X) + \penpar \big(\|X_{-}\|_{\ell_p}^{q} + \|X^\trs
        X \!-\! I_\col\|_{\ell_p}^{\frac{1}{2}}\big)
    \mathrel{:} X \!\in \sset\Big\}.
    \end{split}
\end{equation}

\end{theorem}

Theorem~\ref{th:lexactposorth} presents the exactness of problem~\eqref{eq:nocpen} regarding local minimizers
when~$F$ is locally Lipschitz continuous on~$\sset$, meaning that for any~$\bar{X} \in
\sset$ there exists a constant~$L\in(0,\infty)$ such that~\eqref{eq:lip} holds for all~$X$ and~$Y$
in a certain neighborhood of~$\bar{X}$ in~$\sset$. We will refer to
this~$L$ as a Lipschitz constant of~$F$ around~$\bar{X}$.
The proof of Theorem~\ref{th:lexactposorth} is also given in Appendix~\ref{app:proof}.
\begin{theorem}[Exact penalty~(\ref{eq:nocpen}) with~$F$ being locally Lipschitz continuous]
    \label{th:lexactposorth}
    Let~$\sset\subset\RR^{n\times\col}$ be a set containing~$\stf{n,\col}_+$,~$F\mathrel{:}\sset\to\RR$ be
    a locally Lipschitz continuous function, and~$p\ge 1$ be a constant.
    Suppose that~$0< q_1 \le 1/2$ and~$0 <q_2\le 1/2$.
    For any local minimizer~$X^*$ of~$F$ on~$\stf{n,\col}_+$,~$X^*$ is also a local minimizer of
    \begin{equation}
        \label{eq:nocpenposorth}
        \min \left\{ F(X) +  \penpar(\|X_{-}\|_{\ell_p}^{q_1} +\|X^\trs X - I_\col\|_{\ell_p} ^{q_2})
        \mathrel{:} X\in\sset \right\}
    \end{equation}
    for all~$\penpar > 4L^*\sqrt{\col}\max\Big\{1, \, (n\col)^{\frac{q_1(p-2)} {2p}}\!, \,
    \col^{\frac{q_2(p-2)}{p}}\Big\}$, where~$L^*$ is a Lipschitz constant of~$F$ around~$X^*$.
    Conversely, if~$X^*$ lies in~$\stf{n,\col}_+$ and
    there exists a constant~$\penpar$ such
    that~$X^*$ is a local minimizer of~\eqref{eq:nocpenposorth},
    then~$X^*$ is also a local minimizer of~$F$~on~$\stf{n,\col}_+$.%
\end{theorem}

Suppose that~$p \le 2$. It is noteworthy that the thresholds for~$\penpar$ in
Theorems~\ref{th:exactposorth} and~\ref{th:lexactposorth} are independent of~$n$~(even the dependence on~$\col$
is mild). This is favorable in practice, as~$\col$ can be much smaller than~$n$ in applications.
We also note that the second part of Theorem~\ref{th:lexactposorth} requires~$X^*\in\stf{n,\col}_+$.
This is indispensable without additional assumptions on the problem
structure~(see~\cite[Remark~9.1.1]{Cui_Pang_2021}).

\subsection{The exponents in the penalty term}\label{ssec:exponent}

When~$1 < \col <n$, the requirements on the exponents of~$\|X_{-}\|_\fro$ and~$\|X^\trs X - I_\col\|_\fro$
in Theorems~\mbox{\ref{th:exactposorth} and~\ref{th:lexactposorth}} cannot be relaxed.
This is elaborated in Proposition~\ref{prop:nonexact}, with~$\sset = \RR^{n\times \col}$ being an
example. Similar results can be proved for~$\sset = \stf{n,\col}$ and~$\sset = \RR^{n\times \col}_+$.
\begin{proposition}
    \label{prop:nonexact}
    Suppose that~$1 < \col < n$, $p \ge 1$, $q_1 > 0$, and~$q_2 > 0$.
    Define the function
    $\res(X) = \|X_{-}\|_{\ell_p}^{q_1} + \|X^\trs X - I_\col\|_{\ell_p}^{q_2}$
    for~$X\in\RR^{n\times \col}$.
    There exists a Lipschitz continuous function~$F\mathrel{:} \RR^{n\times \col} \to \RR$ such
    that the following statements hold.
    \begin{enumerate}
    \item\label{it:min}
        $\Argmin\{F(X) \mathrel{:} X\in\stf{n,\col}_+\} =
        \stf{n,\col}_+$.
    \item\label{it:argmin}
        If~$q_1 > 1/2$ or~$q_2 \neq 1/2$, then any~$X^*\in\stf{n,\col}_+$ is not
        a global minimizer of~$F + \penpar \res$
        on~$\RR^{n\times \col}$ for any~$\penpar >0$.
    \item\label{it:arglocmin}If~$q_1 > 1/2$ or~$q_2 > 1/2$, then there exists an~$X^*\in\stf{n,\col}_+$ that is
        not a local minimizer of~$F + \penpar\res$
        on~$\RR^{n\times \col}$ for any~$\penpar >0$.
\end{enumerate}
\end{proposition}
\begin{proof}
    Define
    \begin{equation}
        \nonumber
        F(X) \;=\; -\dist(X,\,\stf{n,\col}_+) \quad \text{for} \quad X\in\RR^{n\times \col}.
    \end{equation}
    Then~$F$ is Lipschitz continuous on~$\RR^{n\times\col}$. We will
    justify~\ref{it:min}--\ref{it:arglocmin} one by one.

    \ref{it:min}~This holds because~$F$ takes a constant value~$0$ on~$\stf{n,\col}_+$.

    \ref{it:argmin}~Assume for contradiction that there exists an~$X^*\in\stf{n,\col}_+$ such
    that~$X^*$ is a global minimizer of~$F + \penpar^* \res$ on~$\RR^{n\times \col}$
    for a certain~$\penpar^* >0$. Then
    \begin{equation}
        \nonumber
        F(X) + \penpar^* \res(X) \;\ge\; F(X^*) + \penpar^*\res(X^*) \;=\; 0
        \quad \text{for all} \quad X \in\RR^{n\times \col}.
    \end{equation}
    By the definition of~$F$, we then have~$\dist(X,\,\stf{n,\col}_+) \le \penpar^*\res(X)$ for
    all~$X \in\RR^{n\times \col}$. Hence~$\res$ defines a global error bound for~$\stf{n,\col}_+$
    relative to~$\RR^{n\times \col}$, contradicting~\ref{it:fro} of Theorem~\ref{th:power}~(note
    that~$\|\cdot\|_{\ell_p}$ and~$\|\cdot\|_\fro$ are equivalent norms).

    \ref{it:arglocmin}~According to~\ref{it:fro} of Theorem~\ref{th:power}, the function~$\res$
    does not define a local error bound for~$\stf{n,\col}_+$ relative to~$\RR^{n\times \col}$.
    Thus there is a sequence~$\{X_k\}\subset \RR^{n\times \col}$ such that%
    \begin{gather}
        \label{eq:ressmall}
        \|(X_k)_{-}\|_\fro + \|X_k^\trs X_k - I_\col\|_\fro \;\le\; k^{-1}, \\
        \label{eq:distbig} \dist(X_k,\,\stf{n,\col}_+) \;>\; k\res(X_k)
    \end{gather}
    for each~$k\ge 1$. According to~\eqref{eq:ressmall}, $\|X_k^\trs X_k\|_\fro \le \sqrt{\col} + k^{-1}$.
    Thus~$\{X_k\}$ has a subsequence~$\{X_{k_{\ell}}\}$ that converges to a certain point~$X^*$.
    Using~\eqref{eq:ressmall} again, we have~$\|X_{-}^*\|_\fro + \|(X^*)^\trs X^* -I_\col\|_\fro
    = 0$, and hence~$X^*\in\stf{n,\col}_+$. It remains to show that~$X^*$ is not a local
    minimizer of~$F+\penpar \res$ for any~$\penpar >0$. Assume for contradiction that~$X^*$ is such
    a local minimizer for a certain~$\penpar^* > 0$.
    Then for all sufficiently large~$\ell$,%
    \begin{equation}
        \nonumber
        F(X_{k_\ell}) + \penpar^* \res(X_{k_\ell}) \;\ge\; F(X^*) + \penpar^* \res(X^*) \;=\; 0.
    \end{equation}
    By the definition of~$F$, we then have~$\dist(X_{k_\ell},\,\stf{n,\col}_+) \le \penpar^* \res(X_{k_\ell})$,
    contradicting~\eqref{eq:distbig}. The proof is complete.
\end{proof}

When~
$\col = 1$ or~$\col = n$,
since the exponents of~$\|X_{-}\|_\fro$ and~$\|X^\trs X - I_\col\|_\fro$
in the error bounds can be increased from~$1/2$ to~$1$, their exponents in the penalty term
of~\eqref{eq:nocpen} can be taken from a larger range while keeping the exactness of~\eqref{eq:nocpen}.
This is briefly summarized in Remark~\ref{rem:r1n}.
\begin{remark}
    \label{rem:r1n}
    Suppose that~$\col = 1$ or~$\col = n$.
    If~$F$ is Lipschitz continuous on~$\sset$, then we can establish a result similar
    to~Theorem~\ref{th:exactposorth} for~$0<q_1 \le 1$ and~$1/2 \le q_2 \le 1$ based on the error
    bound~\eqref{eq:res2bd}.
    When~$F$ is only locally Lipschitz continuous, similar to Theorem~\ref{th:lexactposorth}, the
    exactness of problem~\eqref{eq:nocpen} regarding local minimizers can be established
    if~$0< q_1 \le 1$ and~$0< q_2 \le 1$.
    Proposition~\ref{prop:nonexact} can also be adapted to the case of~$\col = 1$ or~$\col = n$.
    It is also worth noting that~$\stf{n,n}_+$ is precisely the
    set of~$n\times n$ permutation matrices, and
    hence~$\min\{F(X) \mathrel{:} X \in \stf{n,n} _+\}$
    represents optimization problems over permutation matrices.
\end{remark}

\section{Penalty methods for~(\ref{eq:optnoc})}
\label{sec:test}

When~$\stf{n,\col}_\tsign= \stf{n,\col}_+$, problem~(\ref{eq:optnoc}) reduces to the nonnegative orthogonal constrained optimization problem
\begin{equation}\label{L1X-}
\min_{X\in  \stf{n,\col}_{+}} F(X).
\end{equation}
Many papers use penalty methods  for problem~(\ref{L1X-})  with  penalty functions
$\|.\|_\fro^2$, $\|.\|_\fro$ or~$\|.\|_{\ell_1}$ of~$X_{-}$ or~$X^\trs X-I_\col$,
e.g.~\cite{Ahookhosh_etal_2021,Luo_Ding_Huang_Li_2009,Yang_etal_2012,Zass_Shashua_2006}.
However, there is not a satisfactory answer in existing literature whether the penalty problem using
 $\|.\|_\fro^2$, $\|.\|_\fro$ or~$\|.\|_{\ell_1}$  is an exact penalty regarding local and global
 minimizers of problem~(\ref{L1X-}) for a Lipschitz continuous objective function.

In 2024, the authors of~\cite{Qian_Pan_Xiao_2023} proved that the penalty problem
\begin{equation}\label{L1PX-}
\min_{X\in  \stf{n,\col}} F(X) +\mu\|X_{-}\|_{\ell_1}
\end{equation}
is a global exact penalty for problem~(\ref{L1X-})
under the assumption that any global minimizer has no zero rows and~$F$ satisfies a second-order
calmness condition in a neighborhood of any global minimizer of~(\ref{L1X-}).  Moreover, they aimed
to show that such strong assumption  cannot be removed by Example~3.9 in~\cite{Qian_Pan_Xiao_2023}, which is as follows
\begin{equation}\label{IMA2024}
\min_{X\in  \stf{3,2}_+} f(X):=-2X_{11}-2X_{22}-X_{31}-X_{32}.
\end{equation}
The authors of~\cite{Qian_Pan_Xiao_2023} claimed
$X^*=\Bigg[\begin{array}{lll}
1 & 0 &0\\
0 & 1 & 0
\end{array}\Bigg]^\trs$ is a global minimizer of~(\ref{IMA2024}), but is not a solution of the penalty problem
$$
\min_{X\in  \stf{3,2}} f(X) + \mu \|X_{-}\|_{\ell_1}
$$
for any~$\mu>0$.
However,~$X^*$ is not a global minimizer of~(\ref{IMA2024}), as~$f(X^*)=-4$
$>f(\hat{X})= -\sqrt{5}-2$, where~$\hat{X}=\Bigg[\begin{array}{lll}
2/\sqrt{5} & 0 &1/\sqrt{5}\\
0 & 1 & 0
\end{array}\Bigg]^\trs.$  Thus the claim with this example in~\cite{Qian_Pan_Xiao_2023} is wrong.

In this paper, we give a warning for the penalty problem~(\ref{L1PX-}) in the case where the objective function is only Lipschitz continuous.  From Proposition \ref{prop:nonexact}, we know that
there is a Lipschitz continuous function~$f$ such that any global~(local) minimizer  of
(\ref{L1X-}) is not a global~(local)  minimizer of~(\ref{L1PX-}) for any~$\mu>0$.
On the other hand, from Theorem \ref{th:exactposorth} and Theorem \ref{th:lexactposorth}, we know that
\begin{equation}
    \nonumber
\min_{X\in  \stf{n,\col}} F(X) +\mu\|X_{-}\|^q_{\ell_1}
\end{equation}
is an exact penalty problem for~(\ref{L1X-}) regarding global and local minimizers for
$\mu > 5 L\col^\frac{3}{4}$ and~$q\in (0, 1/2],$ where~$L$ is a  Lipschitz constant of~$F$. Our results
provide theoretical warning and guarantee for penalty methods for problem~(\ref{L1X-}).

Note that Remark \ref{rem:r1n} can been extended to the case~$|\posset|+ |\negset|= 1$ or~$|\posset|+|\negset|= n$.
In particular, the penalty problem
\begin{equation}
    \nonumber
\min_{X\in  \stf{n,\col}} F(X) +\mu\|(S\hp X)_{-}\|_{\ell_1}
\end{equation}
is an exact penalty problem of~(\ref{eq:optnoc}) with~$S_{i,1}=1$ and~$S_{i,j}=0,$ for~$j\neq 1, i=1,\ldots, n$.

Consider the following  sparse trace maximization problem~\cite{SChen}
\begin{equation}\label{trace}
\min_{X\in  \stf{n,\col}} -{\rm tr}(X^\trs A^\trs A  X) +\lambda \|X\|_{\ell_1},
 \end{equation}
where~$A \in \RR^{m\times n}$ is a given matrix.
If $A^\trs A$ is a positive or an irreducible non-negative matrix, then by the Perron-Frobenius
theorem, the largest eigenvalue of~$A^\trs A$ is positive and the corresponding  eigenvector is positive. Hence, for a dense nonnegative data matrix A,  it is interesting to consider
\begin{equation}\label{traceS}
\min_{X\in  \stf{n,\col}_S} -{\rm tr}(X^\trs A^\trs A  X) +\lambda \|T\hp X\|_{\ell_1},
 \end{equation}
 with~$S_{i,1}=1$, $S_{i,j}=0,$ $T_{i,1}=0, T_{i,j}=1,$ for~$j\neq 1, i=1,\ldots, n$.  Since the objective function
 of~(\ref{traceS}) is Lipschitz continuous with Lipschitz constant~$L=2\|A\|_2^2+\col\lambda\sqrt{n}$
 over~$\stf{n,\col}$,
  our results show that
 \begin{equation}\label{traceP}
\min_{X\in  \stf{n,\col}}  -{\rm tr}(X^\trs A^\trs A  X) +\lambda \|T\hp X\|_{\ell_1} +\mu\|(S\hp X)_{-}\|_{\ell_1}
\end{equation}
is an exact penalty problem of~(\ref{traceS}) with~$\mu> 5L\col^{\frac{3}{4}}$.

In~\cite{SChen}, Chen et. al proposed a ManPG~(Manifold Proximal Gradient) algorithm to solve the
following nonsmooth
optimization problem
\begin{equation}
    \nonumber
\min_{X\in  \stf{n,\col}}F(X):= f(X) + h(X),
 \end{equation}
 where~$f$ is smooth,~$\nabla f$ is Lipschitz continuous and~$h$ is nonsmooth,  convex and Lipschitz continuous.
 The objective functions in problem~(\ref{trace}) and problem~(\ref{traceP}) satisfy these conditions.
 Numerical results in~\cite{SChen} show that ManPG outperforms some existing algorithms for  solving
 problem~(\ref{trace}). We compare the two models~(\ref{trace}) and~(\ref{traceP}) for sparse trace
 maximization problem by using the code of~\cite{SChen} downloaded from
 \url{https://github.com/chenshixiang/ManPG}, with the same initial points that are randomly generated by the code.

\subsection{Synthetic simulations}

For given~$m, n$, we randomly generated 20 nonnegative matrices and then normalize the columns by
Matlab functions as follows
$$ A={\rm rand}(m,n), \quad  A={\rm normc}(A).$$
For each randomly generated matrix~$A$, we use ManPG to find an approximate solution~$\hat{X}$
of~(\ref{trace}) and~(\ref{traceP}), respectively.
The reconstructed matrix  and its relative reconstruction error~(RRE) and
percentage of explained variance~(PEV)~\cite{XuZ2015} by using~$\hat{X}$ are defined by
\begin{equation}\label{EER}
\hat{A}=A\hat{X}(\hat{X}^\trs \hat{X})^{-1}\hat{X}^\trs, \quad  {\rm RRE}=\frac{\|A-\hat{A}\|_F}{\|A\|_F}, \quad  {\rm PEV}=
\frac{{\rm tr}(\hat{A}^\trs \hat{A})}{{\rm tr}(A^\trs A)}(\times 100\%).
\end{equation}
In Table 1 and Table 2, we report the average values of RRE and PEV of~$\hat{A}$ by using the
randomly generated 20 nonnegative matrices~$A$ for each~$m$ and~$n$
to compare models~(\ref{trace}) and~(\ref{traceP}) with~$\col=10$.

\begin{table}[htb]
{\small     \centering
            \begin{tabular}{c|c|c|c|c|c|c|c}
            \hline
        \multicolumn{8}{c}          {$m=40, \, n=30 $}\\
        \hline
                      $\lambda, \mu$   & $0.6,150$ & $0.6, 170$ & $0.6,190$ & $0.6,200$ & $1,100$ &1, 110&1, 130\\
        \hline
        model (\ref{trace}) & 0.4029  & 0.4029 & 0.4029  &0.4029 & 0.4046 &0.4046 & 0.4046  \\
        model (\ref{traceP})& 0.3999 & 0.3992& 0.3988& 0.3953 & 0.4029 &0.4008& 0.3981   \\
                              \hline
          \multicolumn{8}{c}          {$ \lambda=0.6, \, \mu=100$}\\
       \hline
                      $m,n$   & $50, 25$ & $50, 50$ & $80,  25$ & $80, 40$ & $80,80$ &200, 25&200, 50\\
        \hline
       model  (\ref{trace}) & 0.3811  & 0.4427  & 0.3846 & 0.4315  & 0.4652 & 0.3860 & 0.4464  \\
       model  (\ref{traceP})& 0.3806 & 0.4409 & 0.3843& 0.4284 & 0.4636 &0.3847 & 0.4451 \\
         \hline
                   \end{tabular}
                         \label{Table1}
                         \caption{Comparison on RRE with different~$(m,n, \lambda,\mu)$ by randomly
                         generated~$A$}
}\end{table}

\vspace{-0.2in}
\begin{table}[htb]
{\small     \centering
            \begin{tabular}{c|c|c|c|c|c|c|c}
            \hline
        \multicolumn{8}{c}          {$m=40, \, n=30$}\\
        \hline
                      $\lambda, \mu$   & $0.6,150$ & $0.6, 170$ & $0.6,190$ & $0.6, 200$ & $1,100$ &1, 110&1, 130\\
        \hline
     model    (\ref{trace}) & 0.8376 &  0.8376  &  0.8376  & 0.8376  & 0.8363 &  0.8363    &  0.8363   \\
     model    (\ref{traceP})& 0.8400 & 0.8406& 0.8410 & 0.8410 & 0.8376 &0.8391& 0.8404 \\
                              \hline
          \multicolumn{8}{c}          {$ \lambda=0.6, \, \mu=100$}\\
        \hline
                      $m, n$   & $50, 25$ & $50, 50$ & $80,25$ & $80,40$ & $80,80$ & $200,25$& $ 200, 50$\\
        \hline
     model    (\ref{trace}) & 0.8547  & 0.8040 & 0.8520 & 0.8138  & 0.7836 & 0.8510 & 0.8007 \\
     model    (\ref{traceP})& 0.8551 & 0.8055 & 0.8523& 0.8164& 0.7850 & 0.8520 & 0.8019 \\
         \hline
                   \end{tabular}
                         \label{Table2}
                         \caption{Comparison on PEV with different~$(m,n, \lambda,\mu)$ by randomly
                         generated~$A$}
}\end{table}

\subsection{Numerical results using Yale face dataset}

The Yale Face dataset contains 165 GIF format gray scale images of 15 individuals with 11 images for each subject,
and one for each different facial expression or configuration. From
\url{http://www.cad.zju.edu.cn/home/dengcai/Data/FaceData.html},  we download the~$165 \times 1024$
facial image matrix~$F_{ace}$.
The~$(15 \times (i-1) +t)$th row of~$F_{ace}$ is the~$t$th image of the~$i$th person,
with~$i = 1,\cdots,15$ and~$t = 1,\cdots,11$.
Each row of~$F_{ace}$ defines a~$32\times 32$ nonnegative matrix.
We use  the first 55 rows of~$F_{ace}$, which include 11 images of each of the first five persons,
to get~55 $32\times32$ nonnegative matrices  and then use
Matlab function \textsf{normc} to normalize each of these matrices.

For each~$32\times 32$  matrix~$A$,  we use ManPG to find an approximate solution~$\hat{X}$
of~(\ref{trace}) and~(\ref{traceP}),   respectively.  We compute
the reconstructed matrix~$\hat{A}$  and its RRE and
PEV by using computed~$\hat{X}$ as~(\ref{EER}).

\begin{figure}
	\centering
	\begin{minipage}{0.256\linewidth}
		\centering
		\includegraphics[width = 0.9\linewidth]{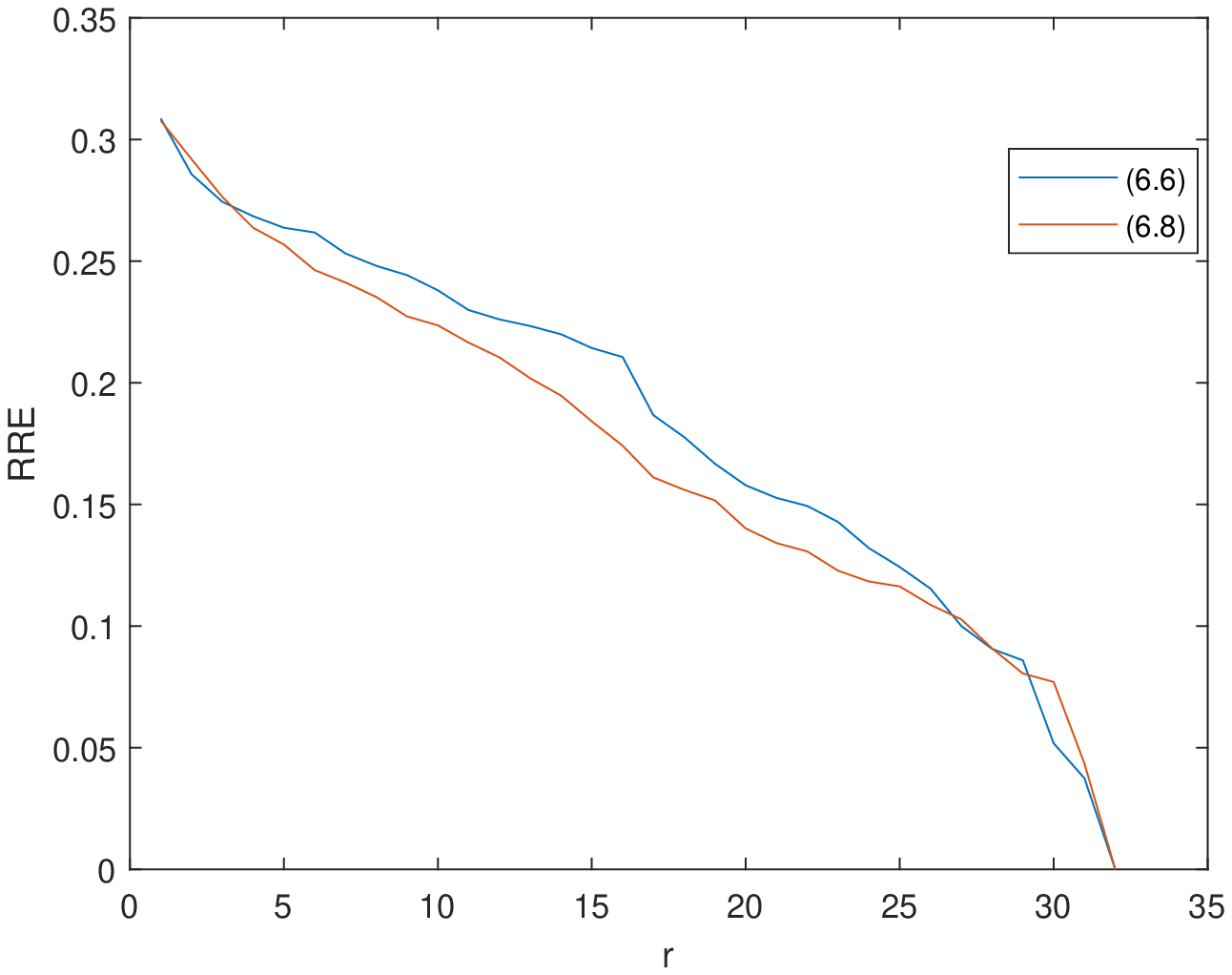}
		\vspace{0.02cm}
		\includegraphics[width = 0.9\linewidth]{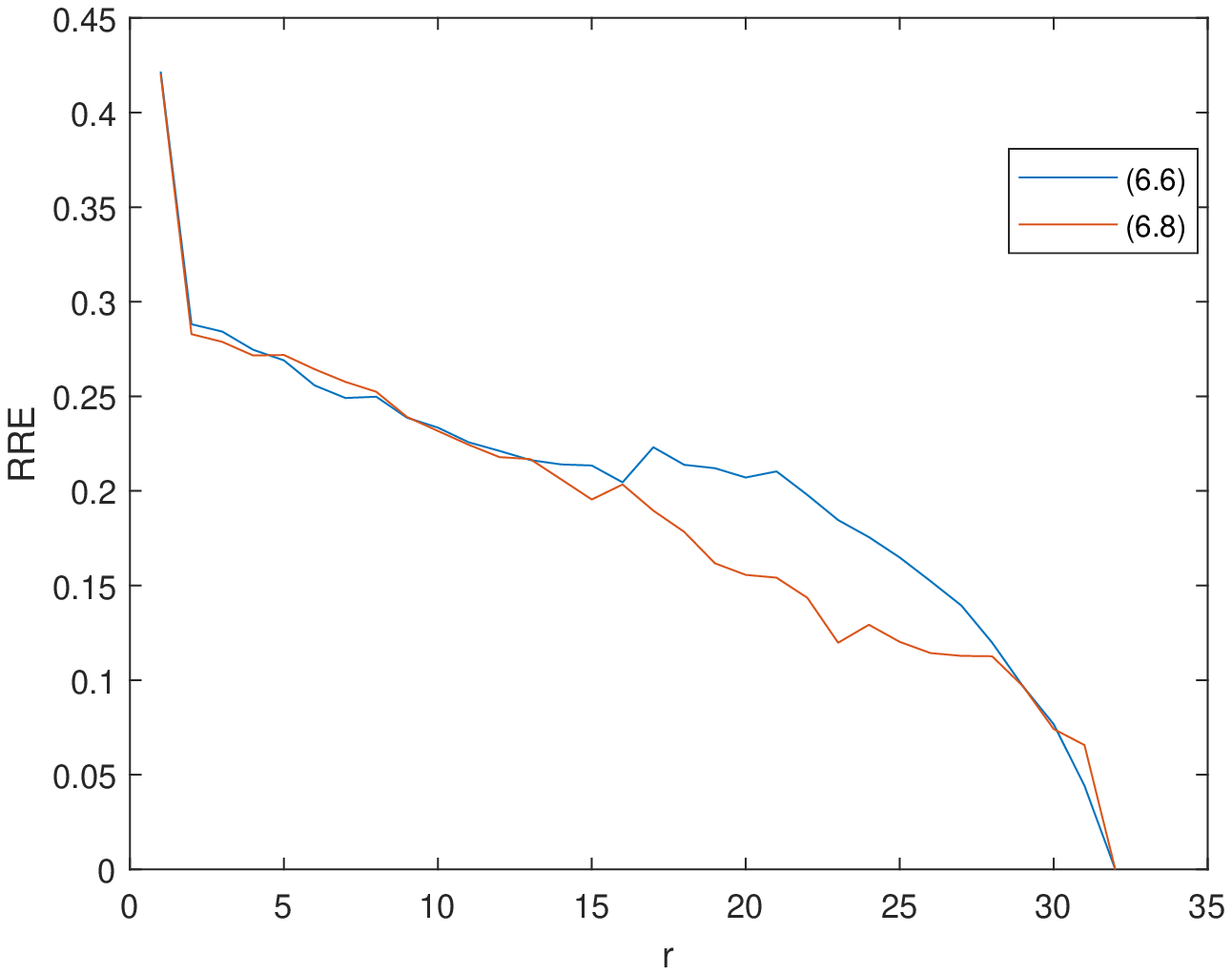}
		\vspace{0.02cm}
		\includegraphics[width = 0.9\linewidth]{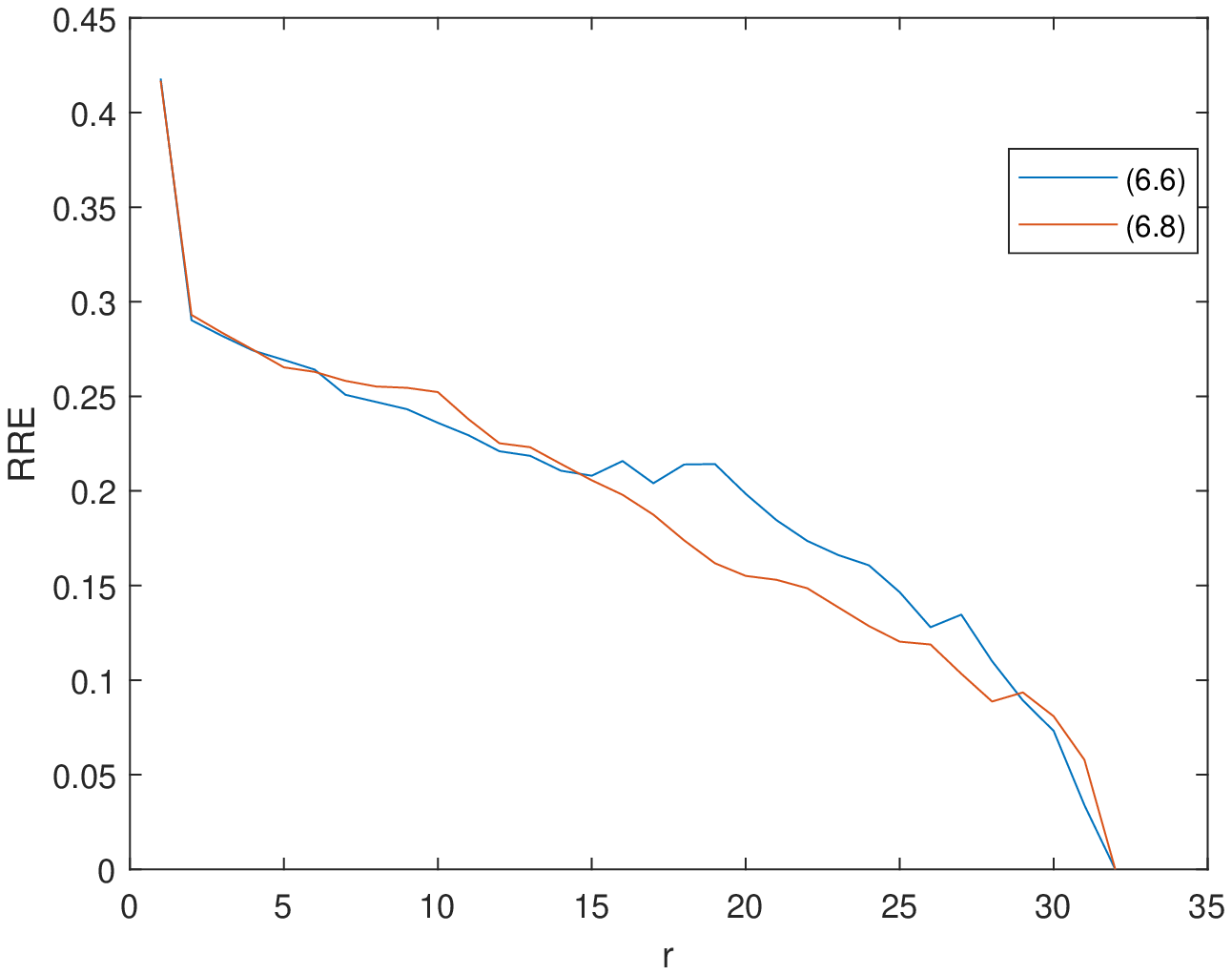}
		\vspace{0.02cm}
		\includegraphics[width = 0.9\linewidth]{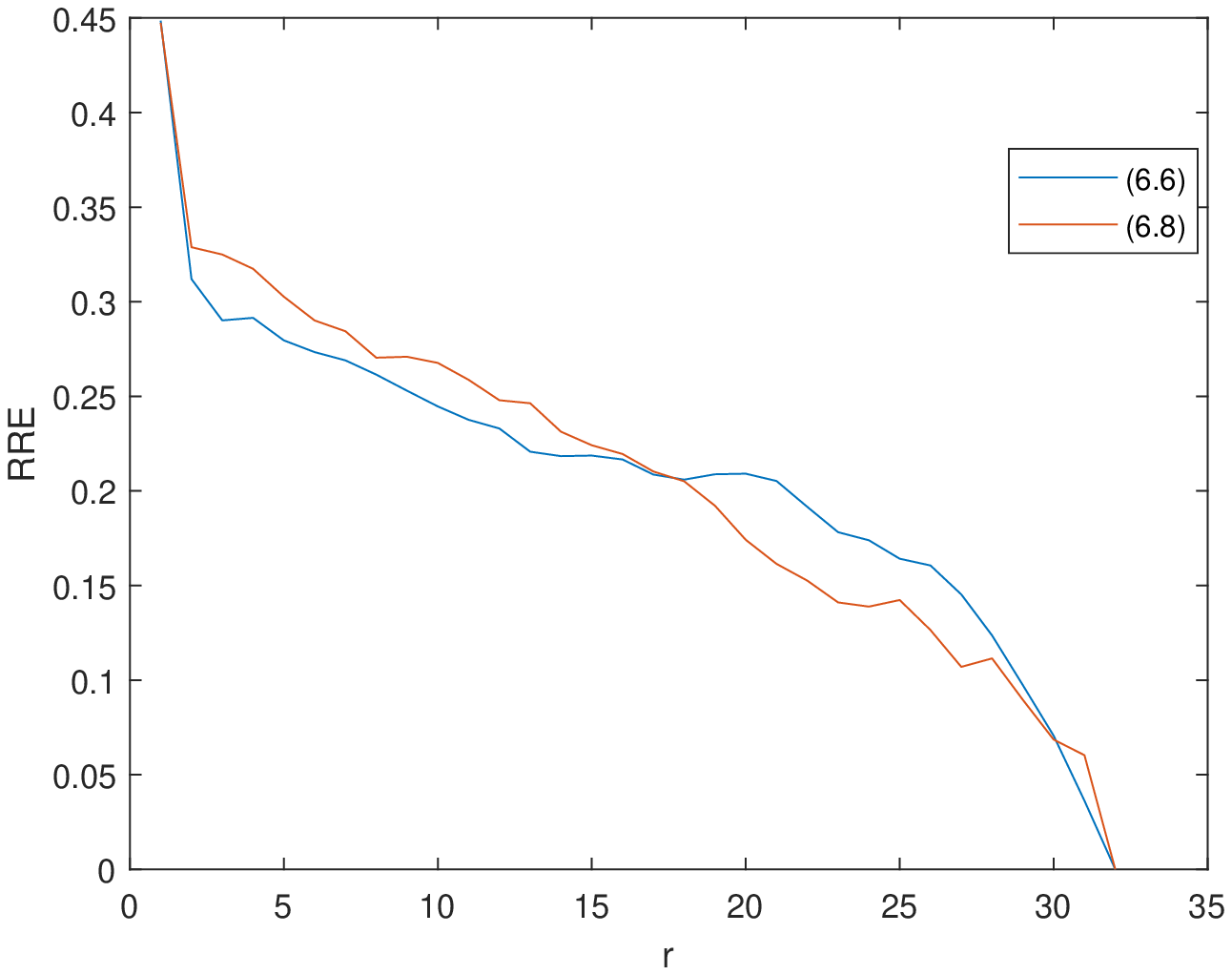}
		\vspace{0.02cm}
		\includegraphics[width = 0.9\linewidth]{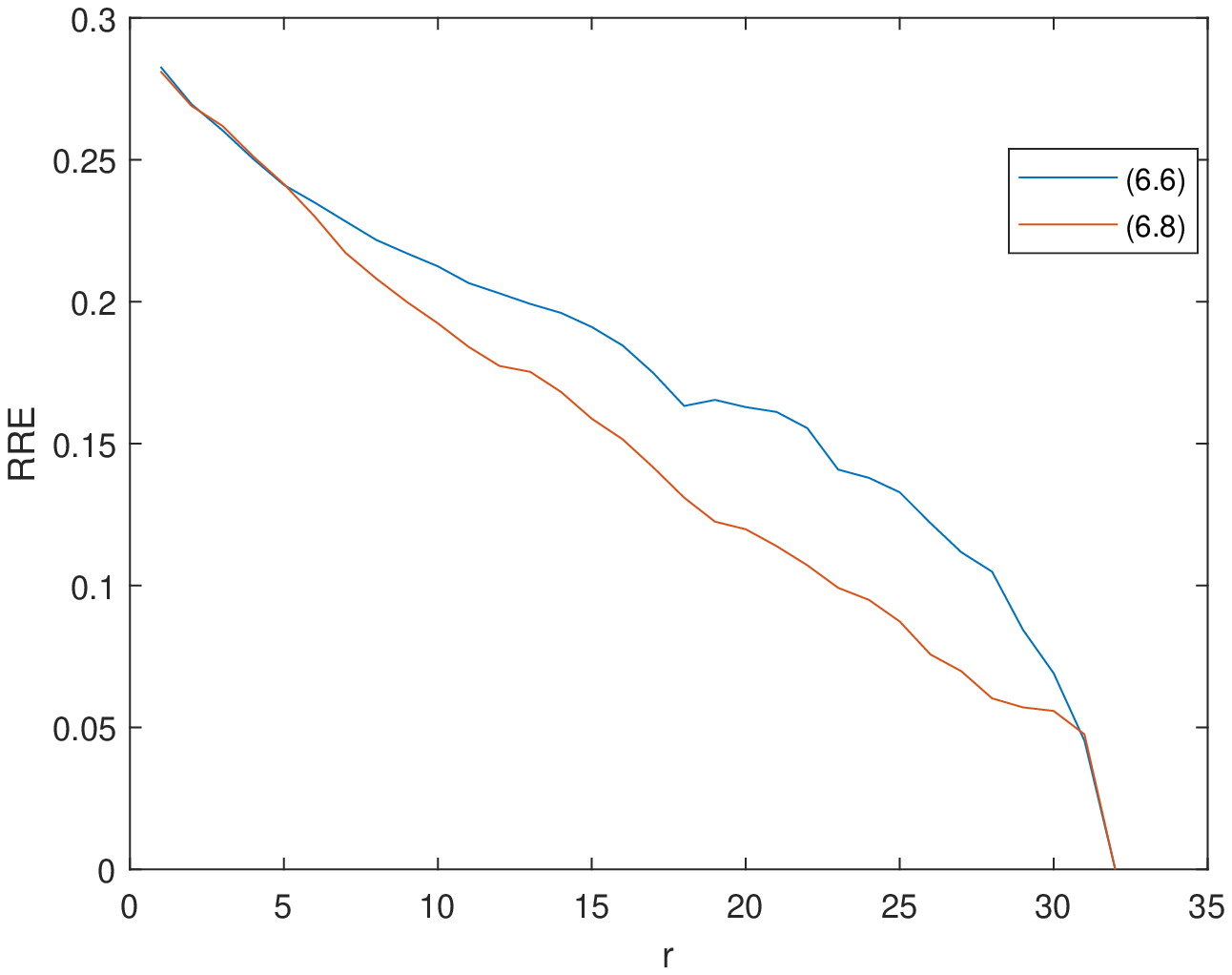}
	\end{minipage}
		\begin{minipage}{0.255\linewidth}
		\centering
		\includegraphics[width = 0.9 \linewidth]{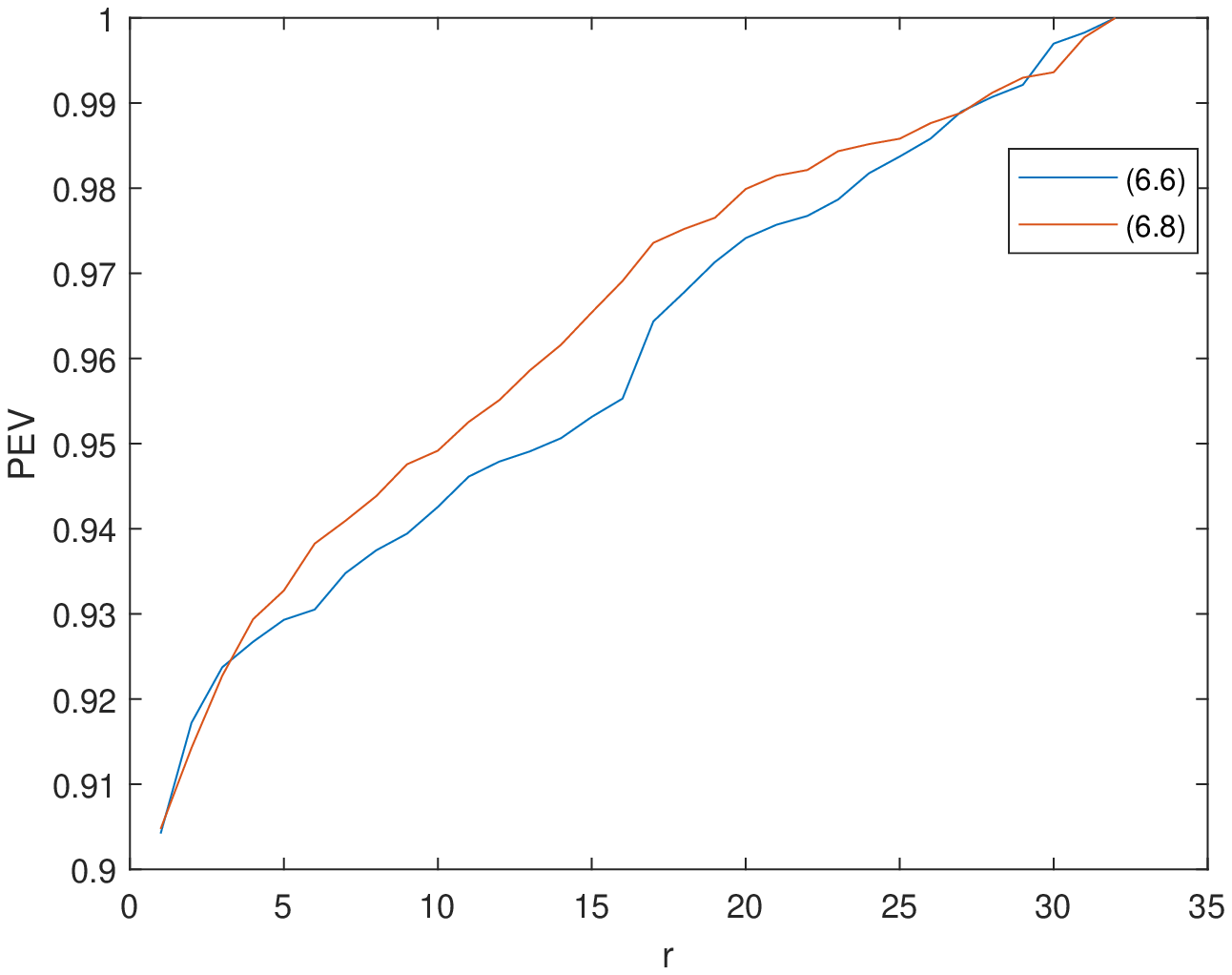}
		\vspace{0.02cm}
		\includegraphics[width = 0.9\linewidth]{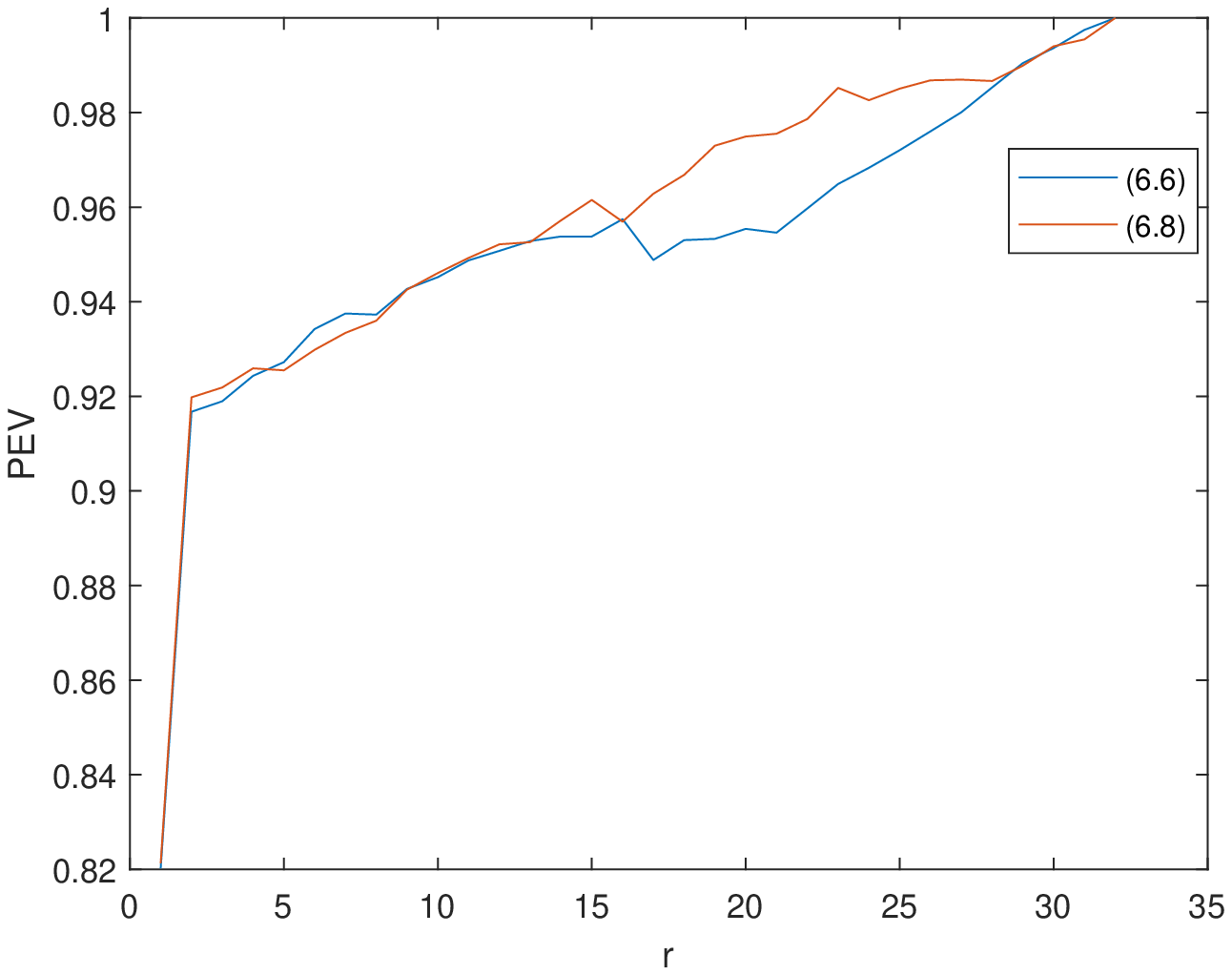}
		\vspace{0.02cm}
		\includegraphics[width = 0.9\linewidth]{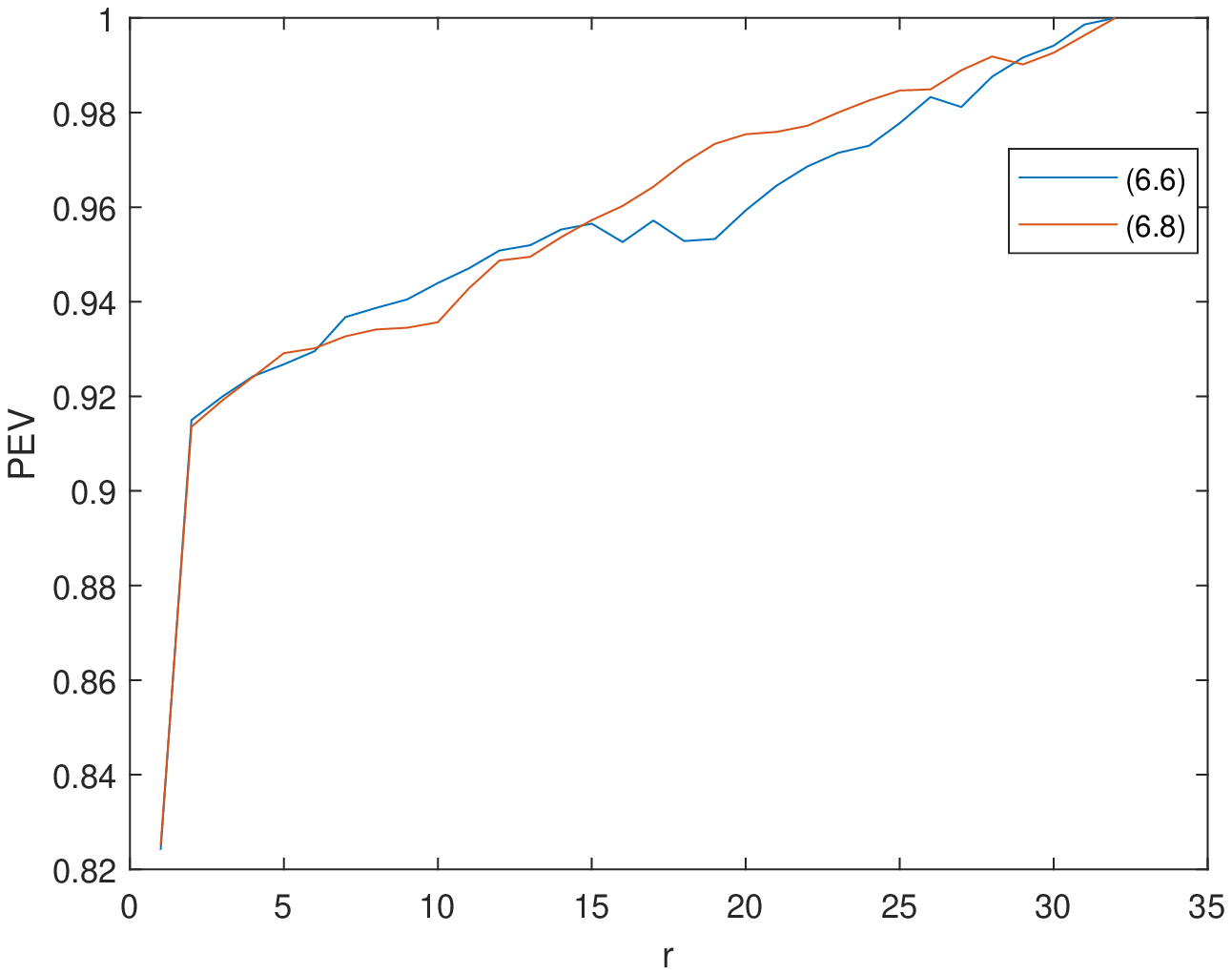}
		\vspace{0.02cm}
		\includegraphics[width = 0.9\linewidth]{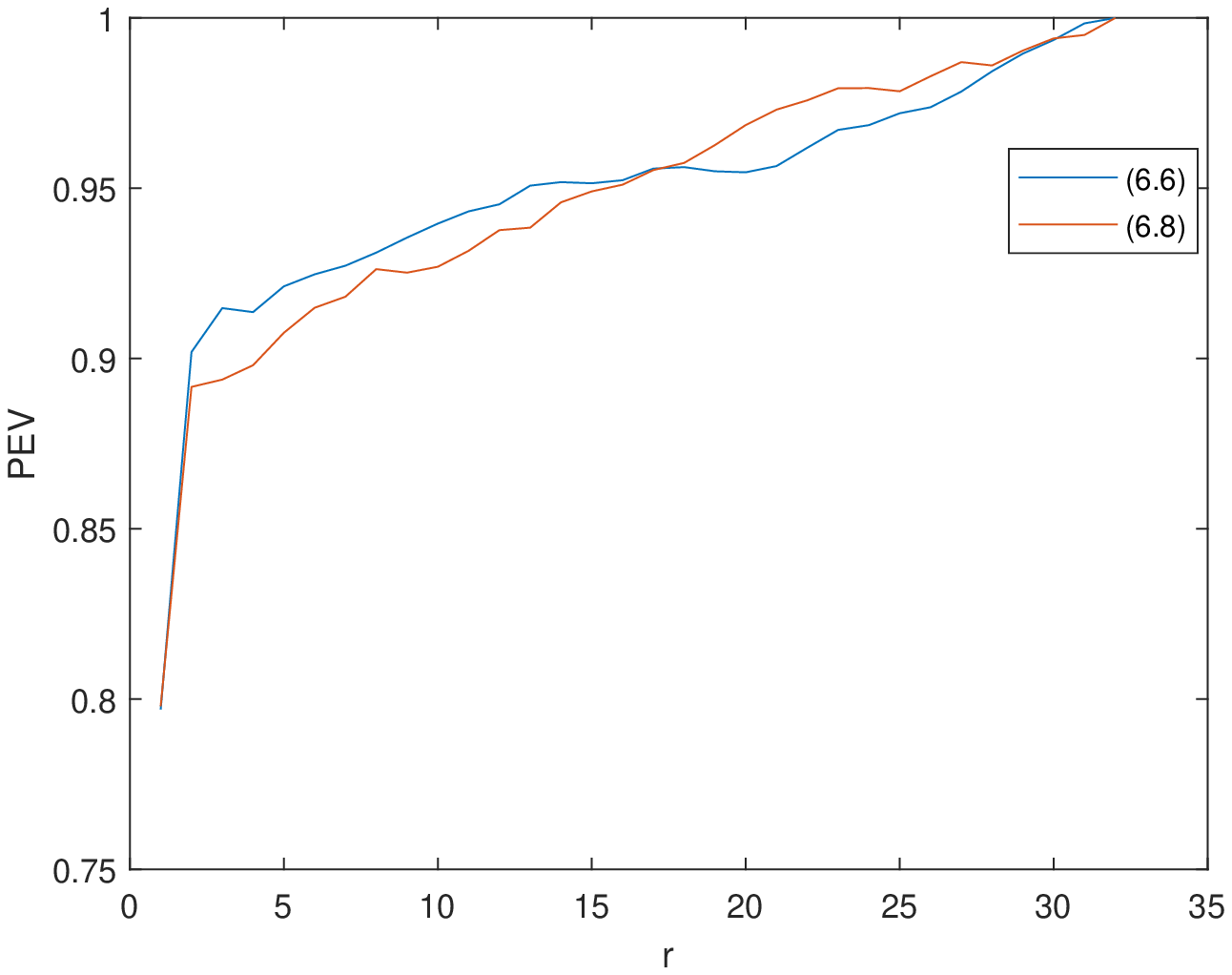}
		\vspace{0.02cm}
		\includegraphics[width = 0.9\linewidth]{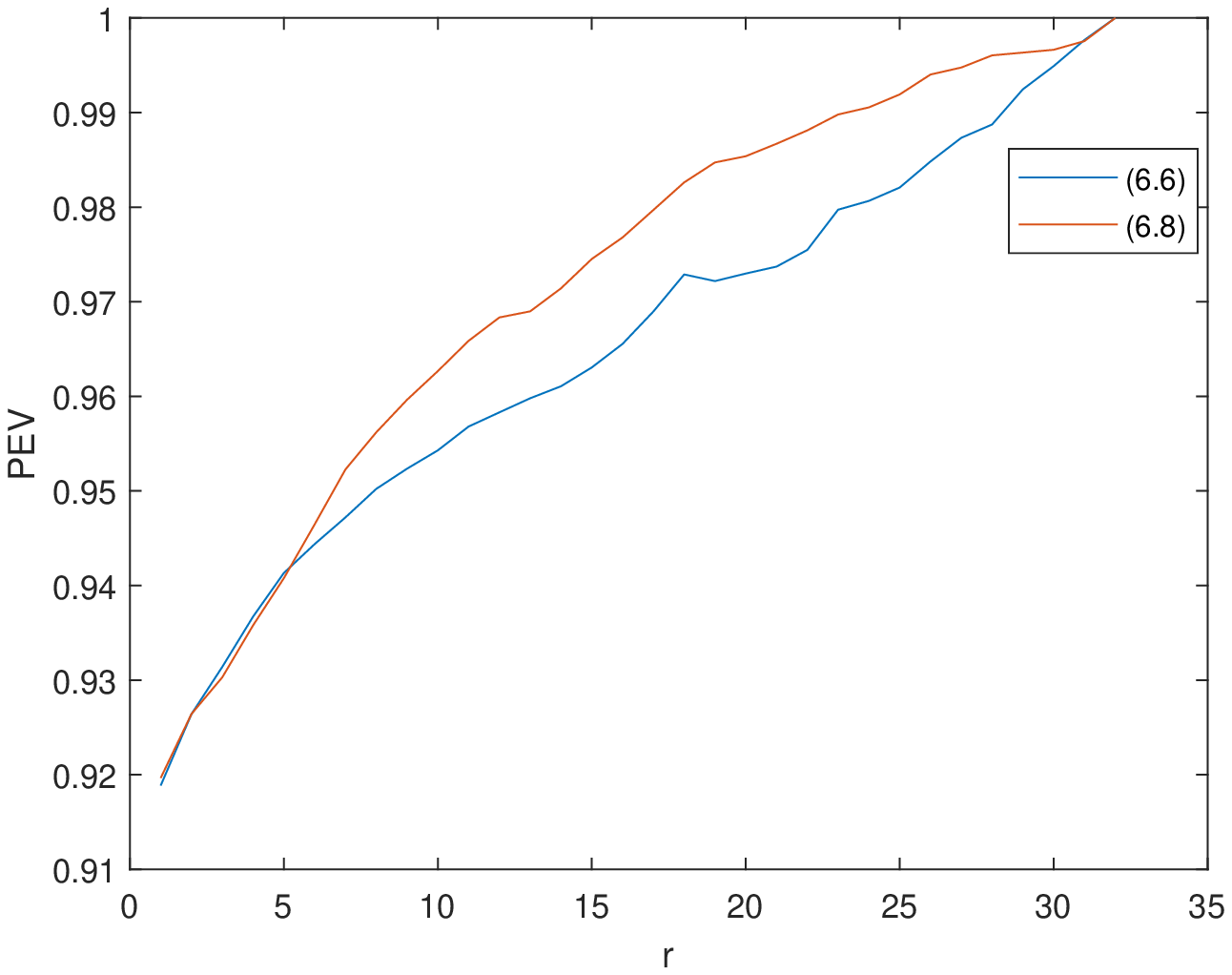}
	\end{minipage}
	\caption{{\small Row~$i$~($i=1,\ldots, 5$)  shows average values of RRE and PEV of the
            reconstructed matrix~$\hat{A}$
using 11 images for the~$i$th person
by models~(\ref{trace}) with~$\lambda=1$ and~(\ref{traceP}) with~$\lambda=1, \mu=6$, respectively,
for~$\col = 1, \cdots, 32$. }}
\end{figure}

From Table 1, Table 2 and Figure 1, we can see that in almost every case, the reconstructed matrix
$\hat{A}$ by model~(\ref{traceP}) has lower values  RRE and higher values PEV than that computed by
model~(\ref{trace}).

\section{Conclusions}
\label{sec:conclusion}

We present the error bounds~(\ref{eq:error1})--(\ref{eq:error3}) with explicit values of~$\nu$ and~$q$ in
Theorems~\ref{th:gerrbdr1},~\ref{th:errbdnr} and~\ref{th:errbd} for~$\stf{n, \col}_\tsign=\stf{n, \col}_+.$
Furthermore,
we show that these error bounds
cannot hold with~$q>{1}/{2}$ when~$1<\col<n$ in Proposition~\ref{prop:tight},
and point out that they cannot hold with~$q>1$ for any~$\col \in \{1, \dots, n\}$ in
Remark~\ref{remark1}.
In Section 4 we present the error bounds~(\ref{eq:error1})--(\ref{eq:error3}) with explicit values of~$\nu$ and~$q$ in
Theorems~\ref{th:errbdsgen1}--\ref{th:errbdsgen} for the sign-constrained Stiefel manifold. The exponent~$q$ in the error bounds is~$1/2$ for any~$\col \in\{1,\dots, n\}$
and can take the value~$1$ for~$|\posset|+|\negset|\in\{1,n\}$.
The new error bounds help us to establish the exactness of penalty problems~(\ref{penalty1})
--(\ref{penalty3}) for problem~(\ref{eq:optnoc}).
Compared with existing results on error bounds for the set~$\stf{n,\col}_+$ and penalty methods for minimization
with nonnegative orthogonality constraints, our results have explicit values of the error bound
parameters and penalty parameters,
and do not need any condition other than the~(local) Lipschitz continuity of the objective function for the exact penalty.
Moreover, exponents in our error bounds are independent of the dimension of the underlying space.

\appendix
\section{Appendix.  Proofs of Theorems~\ref{th:exactposorth} and~\ref{th:lexactposorth}}
\label{app:proof}

We first present the following lemma on a simple inequality between the entry-wise~$\ell_p$-norm and the
Frobenius norm. Its proof is elementary and hence omitted.
\begin{lemma}
    \label{lem:norm}
    For any~$X\in\RR^{n\times\col}$ and any~$p \ge 1$,
\begin{equation}
    \nonumber
    \|X\|_\fro \;\le\; \max\left\{ 1,\; (n\col)^{\frac{p-2}{2p}} \right\}\|X\|_{\ell_p}.
\end{equation}
\end{lemma}

The proofs of Theorems~\ref{th:exactposorth} and~\ref{th:lexactposorth} are as follows.

\begin{proof}[Proof of Theorem~\ref{th:exactposorth}]
    {Define the function~$\pen(X) = \|X_{-}\|_{\ell_p}^{q} + \|X^\trs X -I_\col\|_{\ell_p}^{\frac{1}{2}}$
    for~$X\in\sset$,}
    and set $\nu=5\col^{\frac{3}{4}}\max\Big\{1, \, (n\col)^{\frac{p-2}{4p}}\Big\}$.
        By~\eqref{eq:gerrbdq} and Lemma~\ref{lem:norm}, we have
    \begin{equation*}
        \dist(X,\, \stf{n, \col}_+)
        \;\le\;
        5\col^{\frac{3}{4}} \left(\|X_{-}\|_\fro^{q} + \|X^\trs X -I_\col\|_\fro^{\frac{1}{2}}\right)
        \;\le\; \nu \pen(X)
        \quad  \text{for} \quad X \in \sset.
    \end{equation*}
    For any~$X \in \sset$, setting~$\bar{X}$ to a projection of~$X$ onto~$\stf{n, \col}_+$,
    and combining the \mbox{$L$-Lipschitz} continuity of~$F$ with the above error bound, we have
    \begin{equation}
        \nonumber
        F(\bar{X}) \;\le\;
        F(X) + L\dist(X,\,\stf{n, \col}_+) \;\le\;
        F(X) + \penpar   \pen(X).
    \end{equation}
    This implies that
    \begin{equation}
        \nonumber
        \begin{split}
        \inf\{F(X) \mathrel{:} X \in \stf{n, \col}_+\}
        \;\le\; \inf\{F(X) + \penpar  \pen(X)\mathrel{:}  X\in\sset\}.
        \end{split}
    \end{equation}
    Meanwhile,~$\inf\{F(X) \mathrel{:} X \in \stf{n, \col}_+\}
        \;\ge\; \inf\{F(X) +\penpar  \pen(X)\mathrel{:}  X\in\sset\}$ as~$\pen$ is zero on~$\stf{n, \col}_+\subset\sset$.~Thus%
    \begin{equation}
        \label{eq:infeq1}
        \inf\{F(X) \mathrel{:} X \in \stf{n, \col}_+\}
        \;=\; \inf\{F(X) + \penpar  \pen(X)\mathrel{:}  X\in\sset\}.
    \end{equation}
    For any~$X^*\in   \Argmin\{F(X) \mathrel{:} X \in \stf{n, \col}_+\}$, we have~$\pen(X^*) = 0$ and
    \begin{equation}
        \nonumber
        F(X^*) + \penpar  \pen(X^*) \;=\; F(X^*) \;=\; \inf \{F(X) \mathrel{:} X \in \stf{n, \col}_+\},
    \end{equation}
    which together with~\eqref{eq:infeq1}
    ensures~$X^*\in  \Argmin\{F(X) + \penpar  \pen(X)\mathrel{:} X \in \sset\}$.

    Now take any~$X^* \in \Argmin\{F(X) + \penpar \pen(X)\mathrel{:} X \in
    \sset\}$, and let~$\bar{X}^*$ be a projection of~$X^*$ onto~$\stf{n, \col}_+$. Then we have
    \begin{equation}
        \nonumber
        F(X^*) + \penpar   \pen(X^*)  \;\le\; F(\bar{X}^*)  + \penpar \pen(\bar{X}^*)
        \;=\; F(\bar{X}^*) \;\le\; F(X^*) + \nu L \pen(X^*).
    \end{equation}
   This leads to~$\pen(X^*)=0$,
   as~$\penpar >  \nu L$ and~$\pen(X^*) \ge 0$.  Hence~$X^*$ lies in~$\stf{n, \col}_+$, and
    \begin{equation}
        \nonumber
        F(X^*) \;=\; F(X^*) + \penpar \pen(X^*) \;=\; \inf\{F(X) +\penpar \pen(X)\mathrel{:}  x\in\sset\},
    \end{equation}
    which implies that~$X^*\in  \Argmin\{F(X)\mathrel{:} X \in\stf{n, \col}_+\}$ with the help
    of~\eqref{eq:infeq1}.
    We complete the proof.
\end{proof}

\begin{proof}[Proof of Theorem~\ref{th:lexactposorth}]
    Define the function~$\pen(X) = \|X_{-}\|_{\ell_p}^{q_1} + \|X^\trs X -I_\col\|_{\ell_p}^{q_2}$
    for~$X\in\sset$,
     and set~$\nu= 4\sqrt{\col}\max\Big\{1, \, (n\col)^{\frac{q_1(p-2)} {2p}}\!, \, \col^{\frac{q_2(p-2)}{p}}\Big\}$.
    For any~$X\!\in\stf{n, \col}_+$ and any~$Y\!\in \sset$ such
    that~$\|Y - X\|_\fro < 1/(6\sqrt{\col})$, we have
    \begin{equation*}
        \|Y_{-}\|_\fro + \|\sigma(Y) - \ones\|_2
        \;\le\; \|Y-X\|_\fro + \|\sigma(Y) - \sigma(X)\|_2
        \;\le\; 2\|Y-X\|_\fro
        \;<\; \frac{1}{3\sqrt{\col}},
    \end{equation*}
    where the first inequality is because~$X_{-} = 0$ and~$\sigma(X) - \ones = 0$, while the second
    invokes Lemma~\ref{lem:vonNeumann}.
    Hence~\eqref{eq:lerrbdf} and Lemma~\ref{lem:norm} yield
    \begin{equation}
        \nonumber
        \begin{split}
        &\dist(Y,\, \stf{n, \col}_+)\\
        \;\le\; & 4\sqrt{\col} \left(\|Y_{-}\|_\fro^{q_1} + \|Y^\trs Y - I_\col\|_\fro^{q_2} \right)\\
        \;\le\; & 4 \sqrt{\col} \left(\max\left\{1,\,(n\col)^{\frac{q_1(p-2)}{2p}} \right\}\|Y_{-} \|_{\ell_p}^{q_1}
        + \max\left\{1,\,(\col^2)^{\frac{q_2(p-2)}{2p}} \right\} \|Y^\trs Y - I_\col\|_{\ell_p}^{q_2} \right)\\
                \;\le\; &\nu \pen(Y).
        \end{split}
    \end{equation}

    Given a local minimizer~$X^*$ of~$F$ on~$\stf{n, \col}_+$,
    there exists a~$\delta > 0$ such that~$X^*$ is a global minimizer
    of~$F$ on~$\stf{n, \col}_+\cap \BB(X^*, \delta)$
    and~$F$ is~$L^*$-Lipschitz continuous in the
    same set.

    It suffices to demonstrate that~$X^*$ is a global minimizer of~$F+\penpar\pen$
    on~$\sset\cap \BB(X^*, \delta/2)$
    for all~$\penpar > \nu L^*$.
    Take any point~$Y\in\sset\cap\BB(X^*,\delta/2)$,
    let~$\bar{Y}$ be a projection of~$Y$ onto~$\stf{n, \col}_+$, and
    note that~$\bar{Y}$ lies in~$\BB(X^*,\delta)$, which is because
    \begin{equation}
        \nonumber
        \|\bar{Y}-X^*\|_\fro \;\le\; \|\bar{Y} -Y\|_\fro + \|Y-X^*\|_\fro \;\le\; \|X^*-Y\|_\fro
        + \|Y-X^*\|_\fro \;<\; \delta.
    \end{equation}
    Then, using the fact that~$\pen(X^*) = 0$, we have
    \begin{equation*}
        \nonumber
        F(X^*) + \penpar \pen(X^*) = F(X^*)  \le F(\bar{Y}) \le F(Y) + L^* \dist(Y,\, \stf{n, \col}_+)
        \le F(Y) + \penpar \pen(Y),
    \end{equation*}
    which is what we desire.

    If~$X^*$ is a local minimizer of~$F+ \penpar \pen$ on~$\sset$, and~$X^*$
    happens to lie in~$\stf{n, \col}_+$, then
    \[
    F(X^*)\;=\;F(X^*) + \penpar \pen(X^*) \;\le\;   F(Y)+\penpar \pen(Y)
\]
    for any~$Y$ that is close to~$X^*$ and
    located in~$\stf{n, \col}_+\subset\sset$.
    Hence~$X^*$ is also a local minimizer of~$F$ on~$\stf{n,\col}_+$.
    We complete the proof.
\end{proof}

\bibliographystyle{plain}
\bibliography{signstiefel}
\end{document}